\documentclass{article}

\usepackage{graphicx} 
\usepackage{natbib}
\usepackage{array,epsfig,rotating}
\usepackage{color}
\usepackage[dvipsnames]{xcolor}

\usepackage{hyperref}
\hypersetup{
  colorlinks,
  linkcolor={red!50!black},
  citecolor={blue},
  urlcolor={blue!80!black}
  }

\usepackage{tabularx}
\usepackage{amsmath}
\usepackage{amssymb}
\usepackage{amsfonts}
\usepackage{amsthm}
\usepackage{float}
\usepackage{verbatim}
\usepackage{caption}
\usepackage{subcaption}
\usepackage{setspace}
\usepackage{bbm}
\usepackage[margin=1in]{geometry}
\usepackage[ruled]{algorithm2e} 

\usepackage{tikz}
\usetikzlibrary{shapes, calc, shapes, arrows, positioning}
\tikzstyle{qedge}=[->,thick,black]
\definecolor{myblue}{rgb}{0.0265,    0.6137,    0.8135}
\definecolor{myyellow}{rgb}{0.9290,    0.6940,    0.1250}
\tikzstyle{neuron}=[draw,circle,minimum size=26pt,inner sep=0pt, fill=black!10]
\tikzstyle{hidden}=[draw,circle,minimum size=26pt,inner sep=0pt, fill=white]
\usetikzlibrary{arrows.meta}
\tikzset{>={Latex[width=3mm,length=2mm]}}
\tikzstyle{arr}=[->, thick, black]
\usetikzlibrary{decorations.text}
\usetikzlibrary{shadows,shadings,shapes.symbols}
\tikzset{
    double color fill/.code 2 args={
        \pgfdeclareverticalshading[%
            tikz@axis@top,tikz@axis@middle,tikz@axis@bottom%
        ]{diagonalfill}{100bp}{%
            color(0bp)=(tikz@axis@bottom);
            color(50bp)=(tikz@axis@bottom);
            color(50bp)=(tikz@axis@middle);
            color(50bp)=(tikz@axis@top);
            color(100bp)=(tikz@axis@top)
        }
        \tikzset{shade, left color=#1, right color=#2, shading=diagonalfill}
    }
}

\usepackage{framed}

\theoremstyle{definition}
\newtheorem{theorem}{Theorem}

\newtheorem{lemma}{Lemma}

\newtheorem{example}{Example}

\newtheorem{problem}[theorem]{Problem}

\def\A{\mathbf A}
\def\tilde{\widetilde}

\newcommand{\aaa}{\boldsymbol \alpha}
\newcommand{\bo}{\boldsymbol }
\newcommand{\nnu}{\boldsymbol\nu}

\newcommand{\cc}{\boldsymbol c}

\newcommand{\Y}{\mathbf Y}

\newcommand{\MP}{\mathbb P}

\newcommand{\pa}{\text{\normalfont{pa}}}

\allowdisplaybreaks

\def\r{\mathbf{r}}
\def\P{\mathbb{P}}
\def\R{\mathbb{R}}
\def\E{\mathbb{E}}
\newcommand{\bm}{\mathbf} 
\DeclareMathOperator*{\argmax}{arg\,max}

\DeclareMathOperator{\Var}{Var}
\DeclareMathOperator{\Geom}{Geom}

\newcommand{\Z}{\mathbb{Z}}

\def\spacingset#1{\renewcommand{\baselinestretch}%
{#1}\small\normalsize}
\spacingset{1}

\title{New directions in algebraic statistics:  \\ 
Three challenges from 2023
}
\author{Yulia Alexandr
\and Miles Bakenhus
\and Mark Curiel
\and Sameer K. Deshpande
\and Elizabeth Gross
\and Yuqi Gu
\and Max Hill 
\and Joseph Johnson
\and Bryson Kagy
\and Vishesh Karwa
\and Jiayi Li
\and Hanbaek Lyu 
\and Sonja Petrovi\'c
\and Jose Israel Rodriguez
}
\AtEndDocument{%
  \pagebreak
\noindent YA: \textsc{Department of Mathematics, University of California, Los Angeles, 520 Portola Plaza,
Los Angeles, CA 90095}
\par\textit{E-mail address}: \texttt{yulia@math.berkeley.edu}
\par\textit{URL}: \url{https://yuliaalexandr.github.io/}

\noindent MB: \textsc{Department of Applied Mathematics, Illinois Institute of Technology, Chicago IL 60616, USA.}
\par\textit{E-mail address}: \texttt{mbakenhus@hawk.iit.edu}
\par\textit{URL}: \url{https://github.com/mbakenhus}

\noindent MC: \textsc{Department of Mathematics, University of Hawai`i at M\={a}noa, 2565 McCarthy Mall,
Honolulu, HI 96822}
\par\textit{E-mail address}: \texttt{curielm@hawaii.edu}
\par \textit{URL:} \url{https://sites.google.com/hawaii.edu/markcuriel/home}

\noindent SKD: \textsc{Department of Statistics, University of Wisconsin--Madison,
1300 University Drive,
Madison, WI 53706}
\par\textit{E-mail address}: \texttt{sameer.deshpande@wisc.edu}
\par \textit{URL:} \url{https://skdeshpande91.github.io}

\noindent EG: \textsc{Department of Mathematics, University of Hawai`i at M\={a}noa, 2565 McCarthy Mall,
Honolulu, HI 96822}
\par\textit{E-mail address}: \texttt{egross@hawaii.edu}
\par \textit{URL:} \url{https://math.hawaii.edu/wordpress/egross/}

\noindent YG: \textsc{Department of Statistics, Columbia University, 515 West 116th Street,
New York, NY 10027}
\par\textit{E-mail address}: \texttt{yuqi.gu@columbia.edu}
\par\textit{URL}: \url{https://stat.columbia.edu/~yuqigu/}

\noindent MH: \textsc{Department of Mathematics, University of California, Riverside, 900 University Ave, Riverside, CA 92521}
\par\textit{E-mail address}: \texttt{max.hill1@ucr.edu}
\par\textit{URL}: \url{https://sites.google.com/view/max-hill/}

\noindent JJ: \textsc{Institutionen f\"or Matematik, KTH, SE-100 44 Stockholm, Sweden}
\par\textit{E-mail address}:
\texttt{josjohn@kth.se}

\noindent BK: \textsc{Department of Mathematics, North Carolina State University, 2311 Stinson Dr. , Raleigh, NC 27607}
\par\textit{E-mail address}: \texttt{bgkagy@ncsu.edu}
\par\textit{URL}: \url{https://brysonkagy.github.io}

\noindent VK: \textsc{Department of Statistics, Operations and Data Science, Temple University, Philadelphia PA 19122, USA.}
\par\textit{E-mail address}: \texttt{vishesh@temple.edu}
\par\textit{URL}: \url{https://www.fox.temple.edu/directory/vishesh-karwa-tuk35269}

\noindent JL: \textsc{Department of Statistics and Data Science, University of California, Los Angeles, 520 Portola Plaza,
Los Angeles, CA 90095}
\par\textit{E-mail address}: \texttt{jiayi.li@g.ucla.edu}
\par\textit{URL}: \url{https://jl2ml.github.io/}

\noindent HL: \textsc{Department of Mathematics, University of Wisconsin-Madison, 480 Lincoln Drive, Madison WI 53706-1388, USA.}
\par\textit{E-mail address}: \texttt{hlyu36@wisc.edu}
\par\textit{URL}: \url{https://hanbaeklyu.com/}

\noindent SP: \textsc{Department of Applied Mathematics, Illinois Institute of Technology, Chicago IL 60616, USA.}
\par\textit{E-mail address}: \texttt{sonja.petrovic@iit.edu}
\par\textit{URL}: \url{https://www.sonjapetrovicstats.com}

\noindent JIR: \textsc{Department of Mathematics, University of Wisconsin-Madison, 480 Lincoln Drive, Madison WI 53706-1388, USA.}
\par\textit{E-mail address}: \texttt{jose@math.wisc.edu}
\par\textit{URL}: \url{https://sites.google.com/wisc.edu/jose/}
}
\usepackage[sectionbib]{chapterbib}
\begin{document}

\maketitle

\begin{abstract}
    In the last quarter of a century, algebraic statistics has established itself as an expanding field which uses 
    multilinear algebra, commutative algebra, computational algebra, geometry, and combinatorics to tackle problems in mathematical statistics. These developments have found applications in a growing number of areas, including  biology, neuroscience, economics, and social sciences. 
    
    Naturally, new connections continue to be made with other areas of mathematics and statistics.     
    This paper outlines three
    such connections: to statistical models used in educational testing, to a classification problem for a family of nonparametric regression models, and to phase transition phenomena under uniform sampling of contingency tables.
    We illustrate the motivating problems,  each of which is for algebraic statistics a  new direction, and demonstrate an enhancement of related methodologies. 
\end{abstract}
\section{Introduction}

We illustrate three new research directions in algebraic statistics which share the following common philosophy: they connect algebraic statistics to applied problems from another research area, re-interpret that problem, and illustrate that this new connection is effective toward solving a family of challenges. The choice of the three sets of problems is not made by ranking, but rather by opportunity: namely, in Fall of 2023, the Institute for Mathematics and Statistics Innovation (IMSI) hosted a long program \emph{Algebraic Statistics and Our Changing World}. The program included two-day working group sessions motivated by a problem presented in the `Questions and Consulting seminar'. This paper illustrates three of the new research directions resulting from these interactions. 
In the spirit of the Oberwolfach \emph{Lectures  on Algebraic Statistics} \citep{drton2008lectures} Chapters $6-7$, each section of this survey is self contained.

In Section~\ref{sec:likelihood-geometry},
Yulia Alexandr, Yuqi Gu, Jiayi Li, and Jose Israel Rodriguez study the likelihood geometry of a statistical model motivated by cognitive diagnosis of latent skills in in educational and psychological measurement.  
A Bless model is a discrete statistical model with latent variables.  Bless is an acronym for a ``binary latent clique star forest.'' Identifiability of these models have been previously established \citep{gu2022blessing, gu2023generic} and reparameterizations have leveraged tools from Algebraic Statistics. The new direction here is to study the likelihood geometry of the models for the statistical inference using maximum likelihood estimation. 
Since the Bless model can serve as a building block for identifiable deep generative models  with multiple latent layers \citep{gu2023bp}, studying the likelihood geometry of the Bless model can pave the way for a deeper understanding of these modern powerful generative models.

In Section~\ref{sec:BART}, Mark Curiel, Sameer Despande, Joe Johnson and Bryson Kagy make progress toward identifying (nearly) equivalent regression trees, which represent piece-wise constant step functions (see Figure~\ref{fig:trivial_tree}). 
The new direction aims to leverage ideas from algebraic statistics and combinatorics to improve the Bayesian Additive Regression Trees \citep[BART;][]{Chipman2010} model for nonparametric regression.

In  Section~\ref{sec:transition}, Miles Bakenhus, Elizabeth Gross, Max Hill, Vishesh Karwa, Hanbaek Lyu, and Sonja Petrovi\'c study phase transitions problems on contingency tables through the lens of algebraic statistics. The motivating problem is the appearance of a sharp phase transition in the estimability of the uniform distribution on the space of tables by the hypergeometric distribution.  The threshold for this phase transition is expressed as a condition on the margins of the table, and has been solved for the two-dimensional case. Interpreting the problem from the point of view of algebraic statistics, we propose the generalization of this phenomenon to multi-way tables and partially solve the problem for the $3$-dimensional case. The two distributions are not only combinatorially interesting, but have statistical relevance: sampling from the hypergeometric distribution is used for exact conditional tests of model fit, while sampling from the uniform is needed for performing conditional volume tests under the multinomial sampling scheme.

\subsection*{Acknowledgements}
This research was performed while the authors were
visiting the Institute for Mathematical and Statistical
Innovation (IMSI), which is supported by the National
Science Foundation (Grant No. DMS-1929348). 
The research of YG is partially supported by National
Science Foundation grant DMS-2210796. EG and MC were supported by the National Science Foundation grant DMS-1945584.
The research by JIR is partially supported by the Alfred P. Sloan Fellowship.  
The
Office of the Vice Chancellor for Research and Graduate Education at UW-Madison
with funding from the Wisconsin Alumni Research Foundation partially supports SKD and JIR.
SP is partially supported by the Simons Foundation Collaboration Grant for Mathematicians \#854770 and DOE/SC award \#1010629.

\section{Likelihood Geometry in a Star-forest Model with Dependent Binary Latent Variables}\label{sec:likelihood-geometry}

This section arose from the Questions and Consulting seminar by Yuqi Gu and subsequent discussions amongst Yulia Alexandr, Jiayi Li, and Jose Israel Rodriguez.

\subsection{Blessed models by parameterization}
Consider the \emph{Binary Latent cliquE Star foreSt} (BLESS) model \citep{gu2022blessing} with the following parametrization. Let $\A\in\{0,1\}^K$ denote the latent random vector, and $\Y\in\{0,1\}^p$ denote the observed random vector. We next describe the distribution of $\A$ and $\Y\mid\A$, respectively, to complete the model specification. Assume the binary latent variables can be arbitrarily dependent on each other with the following saturated parametrization:
\begin{align*}
    \MP(\A = \aaa) = \nu_{\aaa}, \quad \forall\aaa\in\{0,1\}^K,
\end{align*}
where $\sum_{\aaa\in\{0,1\}^K}\nu_{\aaa}=1$.
Assume the observed $Y_1,\ldots,Y_p$ are conditionally independent given the latent $\A$, where each $Y_j$ has exactly one latent parent denoted by $A_{\pa(j)}$ (here $\pa(j)\in\{1,\ldots,K\}$). In other words, the bipartite graph from the latent $\A$ to the observed $\Y$ is a star-forest graph, and it follows from the conditional independence property of the graphical model that $\MP(Y_j\mid\A) = \MP(Y_j\mid A_{\pa(j)})$.
Parameterize the conditional distribution of $Y_j\mid A_{\pa(j)}$ as follows:
\begin{align*}
    \MP(Y_j = 1\mid A_{\pa(j)} = 1) = \theta_{j,+},\quad  
    \MP(Y_j = 1\mid A_{\pa(j)} = 0) = \theta_{j,-}.
\end{align*}
Based on the above assumptions, the marginal distribution of the observed random vector $\Y$ can be written as: for all $\bo y\in\{0,1\}^p$, it holds that
\begin{align*}
\MP(\Y = \bo y) 
&= \sum_{\aaa\in\{0,1\}^K} \MP(\A = \aaa) \prod_{j=1}^p \MP(Y_j\mid \A=\aaa)
\\
&= \sum_{\aaa\in\{0,1\}^K}\nu_{\aaa} \prod_{j=1}^p \MP(Y_j\mid A_{\pa(j)}=\alpha_{\pa(j)}) \\
&= \sum_{\aaa\in\{0,1\}^K}\nu_{\aaa} \prod_{j=1}^p \left[\theta_{j,+}^{\alpha_{\pa(j)}} \theta_{j,-}^{1-\alpha_{\pa(j)}}\right]^{y_j}
\left[\left(1-\theta_{j,+}\right)^{\alpha_{\pa(j)}} \left(1-\theta_{j,-}\right)^{1-\alpha_{\pa(j)}}\right]^{1-y_j}.
\end{align*}
We consider the following inequality constraints on the parameters $\bo\Theta=\left\{\theta_{j,+}, \theta_{j,-}:j\in[p]\right\}$ and $\bo\nu=(\nu_{\aaa}:\aaa\in\{0,1\}^K)$:
\begin{align*}
    \theta_{j,+} &> \theta_{j,-}, \quad \forall j\in[p];\qquad
    \nu_{\aaa} > 0, \quad \forall \aaa\in\{0,1\}^K.
\end{align*}
It is known  that when each latent variable $A_k$ has exactly two observed variables as children, the model parameters $\bo\Theta$ and $\nnu$ are generically identifiable \citep{gu2022blessing}. More specifically, in this case,  $\{\theta_{j,+}, \theta_{j,-}\}$ are identifiable if and only if $\A_k$ is not independent from $(A_1,\ldots,A_{k-1}, A_{k+1}, \ldots,A_K)$ (i.e., parameters $\nnu_{\aaa}$ satisfy certain binomial inequalities).
The arbitrary dependence allowed among the $K$ latent variables makes the BLESS model an expressive modeling tool, and also makes it possible to extend it to identifiable deep generative models with multiple latent layers.

Given $N$ i.i.d. observed vectors $\{\bo y_1,\ldots,\bo y_N\}$ in a sample, the MLE  $(\widehat{\bo\Theta}, \widehat\nnu)$ is defined as the maximizer of the likelihood function $\prod_{i=1}^N \MP(\Y = \bo y_i\mid \bo\Theta,\nnu)$.
The identifiability result stated above has a nice consequence that the maximum likelihood estimator (MLE) of parameters $(\bo\Theta,\nnu)$ is statistically consistent as $N\to\infty$. However, given a finite sample of size $N$, the properties of the MLE is not well understood. Moreover, the computation of MLE is often through the iterative EM algorithm, which is sensitive to parameter initialization.

Relevant references in algebraic statistics include: \cite{allman2019maximum,allman2015tensors,fienberg2009,seigal2018mixtures}.

\subsection{Likelihood geometry of Blessed models}

The \textit{maximum likelihood (ML) degree}~\citep[Chapter 7]{Sullivant-book} of an algebraic statistical model $\mathcal{M}$ is the number of complex critical points of $\ell_u$ on the Zariski closure of $\mathcal{M}$ for generic data $u\in\Delta_{n-1}$. The ML degree measures the algebraic complexity of maximum likelihood~estimation.

\subsubsection{A first example and many states}

\begin{example}\label{ex:quartet-tree}
Consider the Blessed model when $K=2$, given by the graph
\begin{center}
\begin{tikzpicture}
    \filldraw (0,1) circle (3pt);
    \filldraw (0,-1) circle (3pt);
    \filldraw (1,0) circle (3pt);
    \filldraw (2,0) circle (3pt);
    \filldraw (3,1) circle (3pt);
    \filldraw (3,-1) circle (3pt);
    \draw (1,0) -- (0,1);
    \draw (1,0) -- (0,-1);
    \draw (1,0) -- (2,0);
    \draw (2,0) -- (3,1);
    \draw (2,0) -- (3,-1);
\end{tikzpicture} 
\end{center}
where all random variables are binary. The model has dimension 11 inside $\Delta_{15}$. Its parametrization is given by
$$p_{ij,kl}=v_{00}a_{i0}b_{j0}c_{k0}d_{l0}+v_{10}a_{i1}b_{j1}c_{k0}d_{l0}+v_{01}a_{i0}b_{j0}c_{k1}d_{l1}+v_{11}a_{i1}b_{j1}c_{k1}d_{l1}.$$

The below code computes the implicit description of the model, utilizing bounded-degree Gr\"obner basis computations in Macaulay2.
\begin{leftbar}
\vspace{-11pt}
\begin{verbatim}
R=QQ[v_(0,0)..v_(1,1),a_(0,0)..a_(1,1),b_(0,0)..b_(1,1),c_(0,0)..c_(1,1),
d_(0,0)..d_(1,1), p_(0,0,0,0)..p_(1,1,1,1),MonomialOrder=>Eliminate 20]
probabilities=toList(p_(0,0,0,0)..p_(1,1,1,1))

par = (i,j,k,l)->(
    summ=0;
    for m from 0 to 1 do (
	for n from 0 to 1 do (
    s=v_(m,n)*a_(i,m)*b_(j,m)*c_(k,n)*d_(l,n);
	    summ=summ+s;
    ););
    return summ;)

gs={}; for p in probabilities do gs=append(gs, p-par((baseName p)#1));
G=ideal(gs); time I=ideal(selectInSubring(1,gens gb(G,DegreeLimit=>15)))
codim I
\end{verbatim}
\vspace{-12pt}
\end{leftbar}

Checking the codimension of the resulting ideal, we find that the model is described by sixteen cubics.
Moreover, this model is precisely the 2-mixture of the $4\times 4$ independence model. This can be seen by realizing the parametrization as a product of a $4\times 2$ matrix,  $2\times 2$ matrix, and $2\times 4$ matrix, as follows:

\[
\begin{bmatrix}
p_{00,00} & p_{00,01} & p_{00,10} & p_{00,11}\\
p_{01,00} & p_{01,01} & p_{01,10} & p_{01,11}\\
p_{10,00} & p_{10,01} & p_{10,10} & p_{10,11}\\
p_{11,00} & p_{11,01} & p_{11,10} & p_{11,11}
\end{bmatrix}=
\begin{bmatrix}
a_{00} b_{00} & a_{01} b_{01} \\
a_{00} b_{10} & a_{01} b_{11}\\
a_{10} b_{00} & a_{11} b_{01}\\
a_{10} b_{10} & a_{11} b_{11}
\end{bmatrix}
\begin{bmatrix}
\nu_{00} & \nu_{01}\\
\nu_{10} & \nu_{11}
\end{bmatrix}
\begin{bmatrix}
c_{00} d_{00} & c_{01} d_{01} \\
c_{00} d_{10} & c_{01} d_{11}\\
c_{10} d_{00} & c_{11} d_{01}\\
c_{10} d_{10} & c_{11} d_{11}
\end{bmatrix}^T
\]

Implicitly, the model is known to be described by all sixteen $3\times 3$ minors of the matrix of joint probabilities $(p_{ij,kl})$, which confirms our earlier computations.
Its ML degree is 191; it was computed in~\cite{HRS2014-mldegree-rank-constraints}.

\end{example}
\subsubsection{Implicit equations for general quartet trees and  blessed  cherry orchard models }
We are interested in the generalization of Example \ref{ex:quartet-tree}, where $k$ binary latent variables form a complete graph, and each latent variable is adjacent to two $n$-category observed variables. Let $Y_{i1}$ and $Y_{i2}$ be the observed variables adjacent to the latent variable $H_i$. We will refer to these models as \textit{blessed cherry orchard models}, denoted by $\mathcal{M}_{k,n}$. Let $T=(v_{i_1\ldots i_k})$ denote the $k$-way $2\times \ldots \times 2$ tensor. For each $i\in[k]$, we will let $M_i$ denote the $n^2\times 2$ matrix with the two columns 
\[
M_i^{(j)}=[\mathbb{P}(Y_{i1}=\ell_1|\alpha_i=j)\cdot\mathbb{P}(Y_{i2}=\ell_2|\alpha_i=j)]_{\ell\in[n]^2},\quad \text{one for each value of }j=0,1.\] 
Note that the parametrization of the model $\mathcal{M}_{k,n}$ can be realized as
$$p_{i_1\ldots i_k}=\nu_{i_1\ldots i_k}
\cdot M_{1}^{(i_1)}\otimes M_{2}^{(i_2)}\otimes\cdots \otimes M_{k}^{(i_k)}.$$
Therefore, $\mathcal{M}_{k,n}$ is the model of $k$-way
$n^2\times n^2\times\cdots \times n^2$ tensors with multilinear rank at most $(2,2,\dots, 2)$. 
Its equations are just minors of flattenings, similar to Example~\ref{ex:quartet-tree}. On the other hand, mixtures of two independence models correspond to border rank at most two matrices, and their equations are much more subtle. Invariants for these cases can be found in \cite{landsberg2012tensors} and \cite{raicu}.

Similarly, when the observed variables have a different number of states, 
the resulting models are multilinear rank $(2,2,\dots,2)$ tensors of different size.
A recursive formula can be derived for these models by the results in 
\citet[Section~4]{RW2017-mixture-independence}.

\subsubsection{EM algorithm for blessed cherry orchard models}
The EM algorithm is the standard method for maximizing the likelihood function on Blessed models. 
Fixed points of the EM algorithm (Algorithm~\ref{alg:em}, also see Algorithm 1 in \cite{gu2022blessing}) on the Blessed cherry orchard model refer to the set of all points $(\bo\Theta^*, \nnu^*)$ where \begin{align*}
    \bo\Theta^{*} &=\arg\max_{\bo\Theta}\mathbb{E}_{\bo A|\bo Y,\bo\Theta^{*},\nnu^{*}}[\log \prod_{i=1}^N \MP(\Y = \bo y_i\mid \bo\Theta^*,\nnu^*), \\
    \nnu^{*} &=\arg\max_{\nnu}\mathbb{E}_{\bo A|\bo Y,\bo\Theta^{*},\nnu^{*}}[\log \prod_{i=1}^N \MP(\Y = \bo y_i\mid \bo\Theta^*,\nnu^*). 
\end{align*} 
Maximizing the log-likelihood on $\mathcal{M}_{k,n}$ is a non-convex optimization problem. 
The output of the EM algorithm $(\bo\Theta^*, \nnu^*)$, with respect to any initialization either lies in the relative interior 
or on the model's boundary.  
If $(\bo\Theta^*, \nnu^*)$ is in the relative interior then $(\bo\Theta^*, \nnu^*)$ is a critical point of the log-likelihood function and the number of such points is counted by the ML degree.
If $(\bo\Theta^*, \nnu^*)$ is on the boundary then $(\bo\Theta^*, \nnu^*)$ is generally not a critical point of the log-likelihood requires studying the ML degree of boundary components like in~\cite{KRS2015-EM}.

\begin{algorithm}[H]
\caption{Expectation-Maximization Algorithm}\label{alg:em}
\KwData{Observation data $\bo Y$, initial parameters $\bo\Theta^{(0)}$, $\nnu^{(0)}$,}
\KwResult{Estimated parameters $\hat{\bo\Theta}$, $\hat{\nnu}$}
\KwIn{Convergence threshold $\epsilon$}
$\bo\Theta \gets \bo\Theta^{(0)}$\;
$\nnu \gets \nnu^{(0)}$\;
$t \gets 0$\;
\Repeat{convergence}{
    $t \gets t + 1$\;
    \textbf{E-step}: Calculate the expected value of the latent variables \;
    $Q(\bo\Theta|\bo\Theta^{(t-1)}) \gets \mathbb{E}_{\bo A|\bo Y,\bo\Theta^{(t-1)},\nnu^{(t-1)}}[\log \prod_{i=1}^N \MP(\Y = \bo y_i\mid \bo\Theta,\nnu)]$\;
    $Q(\nnu|\nnu^{(t-1)}) \gets \mathbb{E}_{\bo A|\bo Y,\bo\Theta^{(t)},\nnu^{(t-1)}}[\log \prod_{i=1}^N \MP(\Y = \bo y_i\mid \bo\Theta,\nnu)]$\;
    
    \textbf{M-step}: Find the parameters that maximize this quantity\;
    $\bo\Theta^{(t)} \gets \arg\max_{\bo\Theta}  Q(\bo\Theta|\bo\Theta^{(t-1)})$\;
    $\nnu^{(t)} \gets \arg\max_{\nnu}  Q(\nnu|\nnu^{(t-1)})$\;
    \textbf{Check for convergence}\;
    \If{$||\bo\Theta^{(t)} - \bo\Theta^{(t-1)}|| < \epsilon$ and $||\nnu^{(t)} - \nnu^{(t-1)}|| < \epsilon$}{
        break\;
    }
}
\end{algorithm}

\subsection{Remaining open questions}

The following open questions about the model and the likelihood geometry are of interest.

First, note that the blessed models we have considered so far have high ML degrees, so there is no closed-form MLE. However, if one imposes symmetries like in \citet[Section 3]{HRS2014-mldegree-rank-constraints} then closed formulas may be derived. Are there any other statistically meaningful restrictions that could guarantee closed-form MLE?
Second, what properties does the EM algorithm have for these models, e.g., do we have the algebraic description of the fixed points of the EM algorithm?
Third, can we perform formal statistical hypothesis tests of the goodness-of-fit of the BLESS model based on the algebraic characterizations? Finding a way to properly sample from logarithmic Voronoi cells~\cite{alexandr-heaton} (as well as describing their boundary) could be useful in this task, as the notion of sufficient statistic is not defined for models with latent variables. Fourth, what is the relationship of blessed models with mixtures of   independence models? Finally, can we determine the boundaries of the image of these models? Are there any nontrivial inequalities?

\newcommand{\skd}[1]{\textcolor{cyan}{\small [skd]: #1}}
\newcommand{\Rlang}{\textsf{R}~}


\definecolor{SkyBlue}{RGB}{14, 118, 188}
\definecolor{BrightRed}{RGB}{223, 82, 78}
\definecolor{Green638}{RGB}{165,255,118} 
\definecolor{myColor1}{RGB}{153,153,153}
\definecolor{myColor2}{RGB}{230,159,0}
\definecolor{myColor3}{RGB}{86,180,233}
\definecolor{myColor4}{RGB}{0,158,115}
\definecolor{myColor5}{RGB}{240,228,66}
\definecolor{myColor6}{RGB}{0,114,178}
\definecolor{myColor7}{RGB}{213,94,0}
\definecolor{myColor8}{RGB}{204,121,167}

\newcommand\numberthis{\addtocounter{equation}{1}\tag{\theequation}} 
\newcommand\numbereqn{\addtocounter{equation}{1}\tag{\theequation}}

\renewcommand{\R}{\mathbb{R}} 
\renewcommand{\E}{\mathbb{E}} 
\def\P{\mathbb{P}} 

\newcommand{\calP}{\mathcal{P}} 
\newcommand{\calQ}{\mathcal{Q}} 
\newcommand{\calF}{\mathcal{F}} 
\newcommand{\calG}{\mathcal{G}}
\newcommand{\calX}{\mathcal{X}}
\newcommand{\calC}{\mathcal{C}}
\newcommand{\calA}{\mathcal{A}}
\newcommand{\calD}{\mathcal{D}}
\newcommand{\calK}{\mathcal{K}}

\newcommand{\cutset}{\calC}
\newcommand{\leftcutset}{\calC_{L}}
\newcommand{\rightcutset}{\calC_{R}}

\newcommand{\ind}[1]{\mathbbm{1}\left( #1 \right)} 
\newcommand{\var}[1]{\textrm{Var}\left( #1 \right)} 
\newcommand{\cov}[2]{\textrm{Cov}\left( #1, #2 \right)} 
\newcommand{\sign}[1]{\textrm{sign}\left(#1\right)} 
\newcommand{\parallelsum}{\mathbin{\|}} 
\newcommand{\kl}[2]{\textrm{KL}\left(#1 \mid \parallelsum \# \right)} 

\newcommand{\normaldist}[2]{\mathcal{N}\left(#1,#2\right)} 
\newcommand{\mvnormaldist}[3]{\mathcal{N}_{#1}\left(#2,#3\right)} 
\newcommand{\gammadist}[2]{\textrm{Gamma}\left(#1,#2\right)} 
\newcommand{\igammadist}[2]{\textrm{Inv.~Gamma}\left(#1,#2\right)} 
\newcommand{\binomialdist}[2]{\textrm{Binomial}\left(#1,#2\right)} 
\newcommand{\berndist}[1]{\textrm{Bernoulli}\left(#1\right)} 
\newcommand{\poisdist}[1]{\textrm{Poisson}\left(#1\right)} 
\newcommand{\hafltdist}[2]{\textrm{half-t}_{\#1}\left(#2\right)} 
\newcommand{\unifdist}[2]{\textrm{Uniform}\left(#1,#2\right)} 
\newcommand{\betadist}[2]{\textrm{Beta}\left(#1,#2\right)} 

\newcommand{\by}{\bm{y}}
\newcommand{\bx}{\bm{x}}
\newcommand{\bz}{\bm{z}}
\newcommand{\bw}{\bm{w}}
\newcommand{\br}{\bm{r}}

\newcommand{\pcont}{p_{\text{cont}}}
\newcommand{\pcat}{p_{\text{cat}}}

\newcommand{\Told}{T_{\text{old}}}
\newcommand{\Tnew}{T_{\text{new}}}

\newcommand{\bY}{\bm{Y}}
\newcommand{\bX}{\bm{X}}

\newcommand{\nx}{\texttt{nx}}
\newcommand{\nxl}{\texttt{nxl}}
\newcommand{\nxr}{\texttt{nxr}}
\newcommand{\nxtil}{\tilde{\nx}}
\newcommand{\nleaf}{\texttt{nleaf}}
\newcommand{\nnog}{\texttt{nnog}}

\newcommand{\bphi}{\boldsymbol{\varphi}}

\newcommand{\Ecal}{\mathcal{E}}

\newcommand{\btheta}{\boldsymbol{\theta}}
\newcommand{\bbeta}{\boldsymbol{\beta}}
\newcommand{\bmu}{\boldsymbol{\mu}}

\newcommand{\ybar}{\overline{y}}
\newcommand{\xbar}{\overline{x}}
\newcommand{\mubar}{\overline{\mu}}

\section{Identifying (nearly) equivalent regression trees}\label{sec:BART}
This section 
grew from the Questions and Consulting seminar by Sameer K. Deshpande and follow up discussions with Mark Curiel, Joseph Johnson, and Bryson Kagy. 
\subsection{Introduction}
\label{sec:introduction}

\textbf{Motivation}. Consider the nonparametric regression problem: given $n$ observations of covariates $\bx \in \R^{p}$ and outcomes $y \in \R$ from the model $y \sim \normaldist{f(\bx)}{\sigma^{2}},$ we would like to estimate the function $f: \R^{p} \to \R.$
Bayesian Additive Regression Trees \citep[BART;][]{Chipman2010} is a Bayesian sum-of-trees model that, at a high-level, approximates $f$ with a large ensemble of regression trees (i.e. piecewise constant step functions).
Usually, $f$ is non-linear and involves complicated high-order interactions.
Such non-linearities and interactions are typically impossible to specify correctly \textit{a priori} using a parametric model. 
Using BART, however, users often obtain extremely accurate predictions of function evaluations along with reasonably well-calibrated uncertainty intervals \textit{without pre-specifying the functional form of $f$ or tuning several hyper-priors}.
The ease-of-use and generally excellent, tuning-free performance have made BART a popular ``off-the-shelf'' tool to be used within larger modeling workflows.

Formally, BART works by simulating draws from a posterior distribution over tree ensembles using Markov chain Monte Carlo.
In each iteration of the sampler, individual trees are grown (by splitting an existing leaf node into two new child nodes) or pruned (by collapsing two leaf nodes to their common parent).
It has been observed empirically --- and recently demonstrated theoretically \citep{Ronen2022_mixing, KimRockova2023_mixing} --- that such local moves result in extremely slow mixing. 
Faster mixing requires more elaborate exploration of the space of regression trees.
In particular, we hope to facilitate faster mixing by directly transitioning between trees that provide (nearly) identical fits to the data.

\subsubsection{Setting \& notation}

To motivate this idea, consider the slightly simpler setting in which we approximate the function $f(\bx)$ with a single binary regression tree and know the residual variance $\sigma^{2}.$
Formally, a \textit{regression tree} is a pair $(T,\bmu)$ consisting of (i) a finite, rooted binary decision tree $T$ containing several terminal or \textit{leaf} nodes and several non-terminal or \textit{decision} nodes and (ii) a collection $\bmu$ of scalars, one for each leaf node in $T.$
Every non-terminal node in $T$ is connected to two children nodes, a left child and a right child.
Further, associated to every non-terminal node is a decision rule of the form $X_{j} < c$ where $c \in (0,1).$

Given a decision tree $T$ and any point $\bx = (x_{1}, \dots, x_{p}) \in [0,1]^{p},$ we can trace a path from the root to a leaf by following the decision rules.
Specifically, starting from the root, whenever the path reaches a decision rule $X_{j} < c,$ it proceeds to the left if $x_{j} < c$ and to the right otherwise.
We will restrict attention only to those decision trees that partition $[0,1]^{p}$ in the sense that (i) every leaf node is reached by the decision-following path of at least one $\bx \in \calX$ and (ii) the decision-following path of every $\bx \in [0,1]^{p}$ reaches a single, unique leaf.
Given $(T, \bmu)$ and a point $\bx \in [0,1]^{p},$ let $\ell(\bx;T)$ denote the leaf reached by $\bx$'s decision-following path.
By associating each leaf of $T$ with its own scalar, the regression tree $(T, \bmu)$ represents a piecewise constant function of $[0,1]^{p}.$
Formally, we introduce the evaluation function $g(\bx; T, \bmu) = \mu_{\ell(\bx)},$ which returns the element of $\bmu$ associated with the leaf reached by $\bx$'s decision-following path.
Additionally, given $n$ points $\bx_{1}, \ldots, \bx_{n} \in [0,1]^{p}$ let $I_{\ell}(T) = \{i: \ell(\bx_{i};T) = \ell\}$ contain the set of indices of the points that reach leaf $\ell.$

With this notation in hand, we can define the single-tree Bayesian model with known residual variance $\sigma^{2}$
\begin{align*}
y_{i} \vert T, \bmu, &\sim \normaldist{g(\bx_{i}; T,\bmu)}{\sigma^{2}} \quad \text{for $i = 1, \ldots, n$} \\
\mu_{\ell} \vert T &\sim \normaldist{0}{\tau^{2}} \text{for $\mu_{\ell} \in \bmu$} \\
T &\sim \Pi(T)
\end{align*}
where $\sigma, \tau, \nu, \lambda > 0$ are fixed positive constants and $\Pi(T)$ is the decision tree prior used in \cite{Chipman2010}.

Under this model, we can compute the marginal likelihood of the decision tree $T$ in closed-form:
$$
p(\by \vert T) \propto \prod_{\ell}{\exp\left\{ \text{stuff depending only on $I_{\ell}(T)$'s \& $y$'s}\right\}}
$$
Importantly, the marginal likelihood of $T$ depends on the decision tree only through the partition of the points $\{\bx_{1}, \ldots, \bx_{n}\}$ it induces.
We say that two decision trees $T$ and $T'$ are \textit{equivalent} if they induce the same partition of $\{\bx_{1}, \ldots, \bx_{n}\}.$

We study the following questions: given a single decision tree $T$ and a collection of points $\bx_{1}, \ldots, \bx_{n} \in \R^{p},$ can we
\begin{enumerate}
\item{enumerate or characterize the equivalence class of trees that induce the exact same partition of the points?}
\item{sample uniformly from the set of trees inducing the same partition?}
\item{enumerate or characterize the set of trees that induce partitions that are close (in some sense) to the one induced by the tree}
\item{sample uniformly from the set of trees inducing nearly the same partition?}
\end{enumerate}

Resolving these questions will enable construction of more efficient MCMC sampling techniques for fitting BART.
Beyond the motivating Bayesian context, however, answering these questions will more generally facilitate uncertainty quantification about tree models.
As one example, suppose we fit a complicated machine learning model (e.g. a deep neural network) to $(\bx, y)$ data to obtain predicted outcomes $\hat{y}.$
Even if the fitted model is difficult to interpret, we can nevertheless obtain a much more interpretable \textit{approximation} by training a regression tree model to the pairs $(\bx, \hat{y}).$
The set of (nearly) equivalent trees provides one avenue to quantify uncertainty about the interpretation of the original fitted model.

\subsection{Preliminaries \& special cases}
\begin{figure}[H]
\centering
\begin{subfigure}{0.35\textwidth}
\centering
\includegraphics[width = \textwidth]{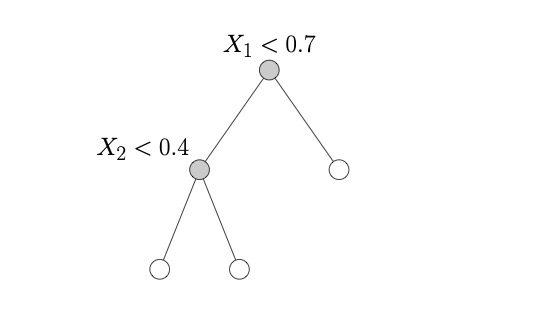}
\end{subfigure}
\begin{subfigure}{0.2\textwidth}
\centering
\includegraphics[width = \textwidth]{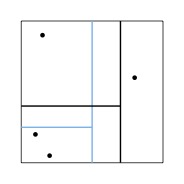}
\end{subfigure}
\begin{subfigure}{0.35\textwidth}
\centering
\includegraphics[width = \textwidth]{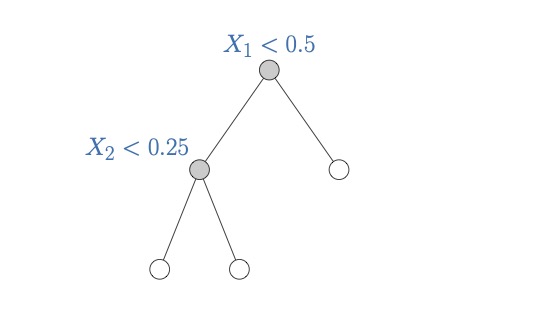}
\end{subfigure}
\caption{We can form trivially equivalent trees by perturbing decision boundaries. Here the tree topology for both trees is the same, but the decision rules for the tree on the left and right are different but they are equivalent because they differentiate the data points in the same way. }
\label{fig:trivial_tree}
\end{figure}

Before proceeding, we verify that equivalent trees do exist generally: given a tree $T$ and set of points $\bx_{1}, \ldots, \bx_{n},$ we can trivially obtain equivalent trees by moving the decision boundaries between the data points; see Figure~\ref{fig:trivial_tree} for an example with $p = 2.$

\begin{figure}[H]
\centering
\begin{subfigure}{0.35\textwidth}
\centering
\includegraphics[width = \textwidth]{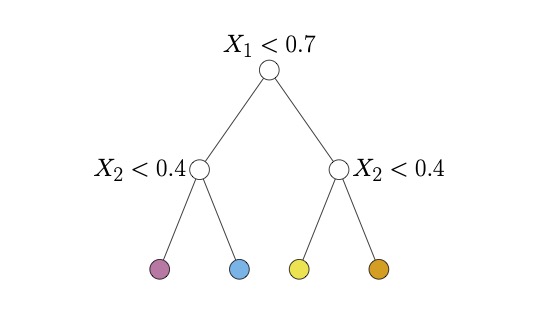}
\end{subfigure}
\begin{subfigure}{0.2\textwidth}
\centering
\includegraphics[width = \textwidth]{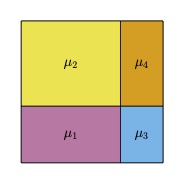}
\end{subfigure}
\begin{subfigure}{0.35\textwidth}
\centering
\includegraphics[width = \textwidth]{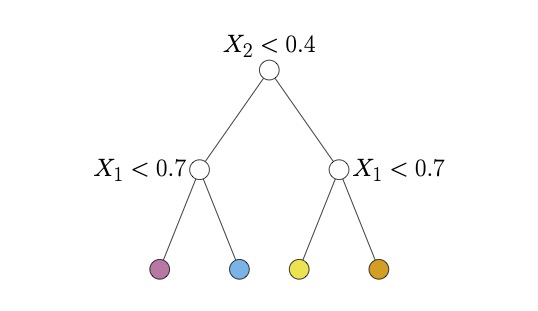}
\end{subfigure}
\caption{For full decision trees with the same rule at each level, we can form equivalent trees by permuting the decision rules across levels.}
\label{fig:full_tree}
\end{figure}

Henceforth, we will focus instead on identifying equivalent trees that do not change the decision boundaries.
To this end, consider first the case where $T$ is (i) a full binary tree of depth $D$ containing $2^{D}$ leaf nodes and (ii) the same decision rule is used at every decision node at depth $d.$
Given such a tree, we can form an equivalent tree by permuting the decision rules across the levels.
Figures~\ref{fig:full_tree} shows an example of equivalent trees of depth $D = 2$ with $p = 2$ and the associated partition of $[0,1]^{2}.$

\subsubsection{The $p = 1$ setting}

\begin{figure}[H]
\centering
\includegraphics[width = 0.7\textwidth]{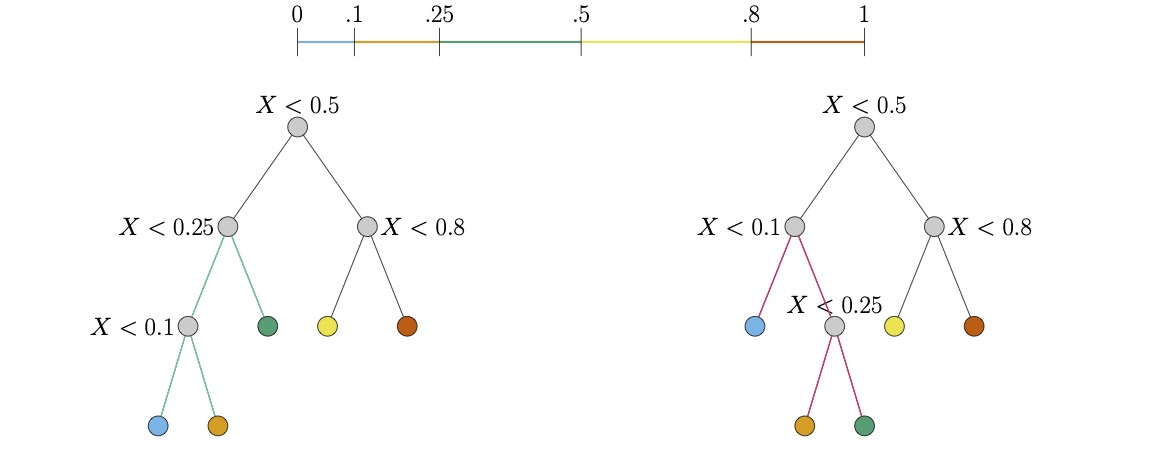}
\caption{We can form equivalent trees by replacing the green sub-tree on the left with the equivalent red sub-tree on the right. Doing so amounts to partitioning $[0,0.5)$ into $[0, 0.1), [0.1, 0.25),$ and $[0.25, 0.5)$ in different orders.}
\label{fig:one_dim_recursion}
\end{figure}

Now suppose that $p = 1$ and that $T$ contains $L \geq 3$ leaf nodes and $L-1$ decision nodes.
In this case, finding equivalent trees is not as simple as permuting decision rules across the levels of $T.$
Notice, however, that we can form an equivalent tree by replacing sub-trees of $T$ with equivalent sub-trees (see Figure~\ref{fig:one_dim_recursion} for an example). 
To see that the two trees are equivalent, notice that original sub-tree partitions partitions the interval $[0,0.5) = [0, 0.25) \cup [0.25, 0.5)$ and then partitions  $[0, 0.25) = [0,0.1) \cup [0.1, 0.25).$
The equivalent sub-tree first partitions $[0, 0.5) = [0, 0.1) \cup [0.1, 0.5)$ and then partitions $[0.1, 0.5) = [0.1, 0.25) \cup [0.25, 0.5).$
Essentially, replacing a sub-tree with an equivalent tree amounts to arranging the decision boundaries in that sub-tree in a different order. 

The example of Figure~\ref{fig:one_dim_recursion} gives us a strategy for counting the number of equivalent trees.
To this end, suppose that $T$ contains $L-1$ internal, decision nodes and let $0 < c_{1} < \cdots c_{L-1} < 1$ be the ordered decision boundaries.
For any collection of $\ell-1$ consecutive decision boundaries let $C_{\ell-1}$ count the number of equivalent trees that contain $\ell-1$ decision nodes associated with those decision boundaries.
Immediately, we know that $\tilde{C}_{\ell-1} < C_{\ell-1}$, the $(\ell-1)^{\text{st}}$ Catalan number, which counts the total number of binary trees with $\ell-1$ internal nodes.

\begin{figure}[H]
\centering
\begin{subfigure}{0.3\textwidth}
\centering
\includegraphics[width = \textwidth]{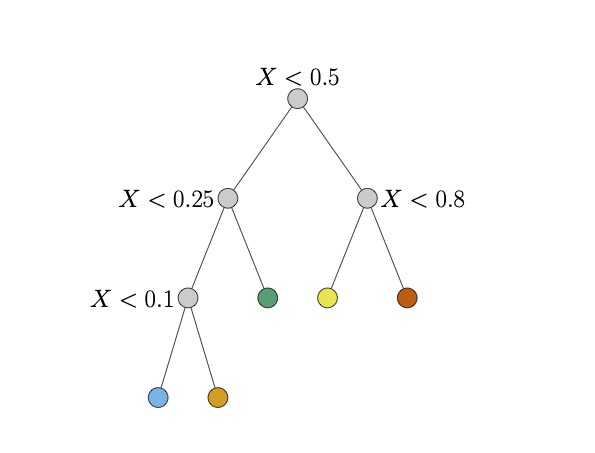}
\caption{}
\label{fig:one_dim_original}
\end{subfigure}
\begin{subfigure}{0.3\textwidth}
\centering
\includegraphics[width = \textwidth]{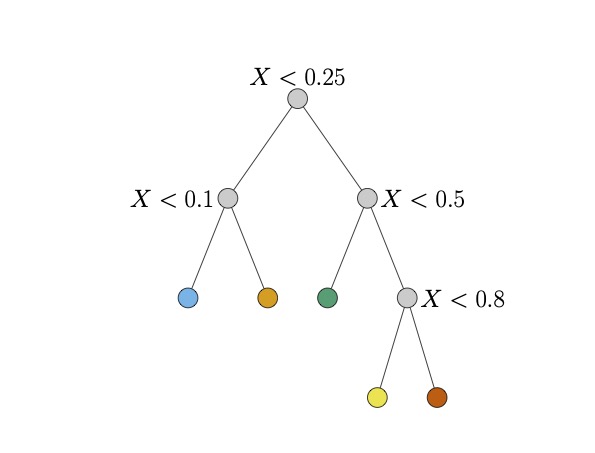}
\caption{}
\label{fig:one_dim_alt1}
\end{subfigure}
\begin{subfigure}{0.3\textwidth}
\centering
\includegraphics[width = \textwidth]{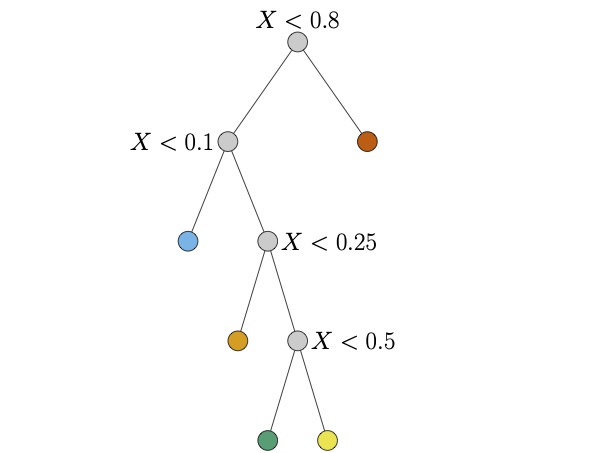}
\caption{}
\label{fig:one_dim_alt2}
\end{subfigure}
\caption{Three trees equivalent to the ones in Figure~\ref{fig:one_dim_recursion}.}
\label{fig:one_dim}
\end{figure}

Further suppose that we have enumerated all $\tilde{C}_{\ell-1}$ such trees for every collection of $\ell - 1$ consecutive boundaries with $\ell < L.$
We can form a new tree $T^{\star}$ as follows:
\begin{enumerate}
\item{Initialize $T^{\star}$ to be just the root node}
\item{Pick one decision boundary $c_{k}$ to associate to $T^{\star}$'s root}
\item{Draw one of the $\tilde{C}_{k-1}$ trees with $k-1$ decision boundaries $\{c_{1}, \ldots, c_{k-1}\}.$ Call it $T^{\star}_{L}$}
\item{Draw one of the $\tilde{C}_{L-k-1}$ trees with $L-k-1$ decision boundaries $\{c_{k+1}, \ldots, c_{L-1}\}.$ Call it $T^{\star}_{R}$}
\item{Connect the roots of $T^{\star}, T^{\star}_{L}$ and $T^{\star}_{R}$ so that $T^{\star}_{L}$'s root is the left child of $T^{\star}$'s root and $T^{\star}_{R}$'s root is the right child of $T^{\star}$'s root.}
\end{enumerate}
Using this process, we have
$$
\tilde{C}_{L-1} \geq \sum_{k = 1}^{L-1}{\tilde{C}_{k-1}\tilde{C}_{L-k-1}}.
$$
A strong induction argument show that, in fact, $\tilde{C}_{L-1} = C_{L-1}$.
So when $p = 1,$ there are a Catalan number of equivalent trees with the same number of leaf nodes.
Further, given any tree with $L$ leaf nodes, we can form every other equivalent tree with the same number of leaf nodes and decision boundaries.
Figure~\ref{fig:one_dim} shows some examples with $L = 5.$

\subsection{The $p = 2$ setting}

\begin{figure}[H]
\centering
\includegraphics[width = 0.6\textwidth]{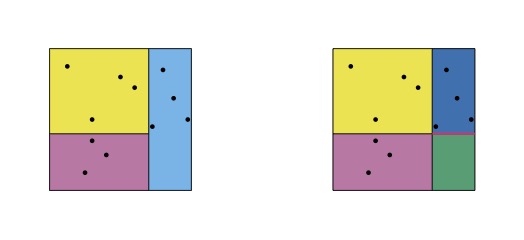}
\caption{Extending the cut along $x_{2} = 0.4$ does not change the underlying partition of the data.}
\label{fig:p2_partition}
\end{figure}

Consider the partitions of $[0,1]^{2}$ and $10$ data points on the left of Figure~\ref{fig:p2_partition}. Notice that we can further partition the blue rectangle $[0.7,1] \times [0,1]$ into $[0.7, 1] \times [0, 0.4) \cup [0.7,1] \times [0.4, 1]$ \emph{without changing the partition of the data points}. That is, we can obtain the partition on the right by ``extending'' the cut separating pink and yellow rectangles across through the blue rectangle; this extended cut is highlighted in the partition on the right of Figure~\ref{fig:p2_partition}.

\begin{figure}[h]
\centering
\includegraphics[width = \textwidth]{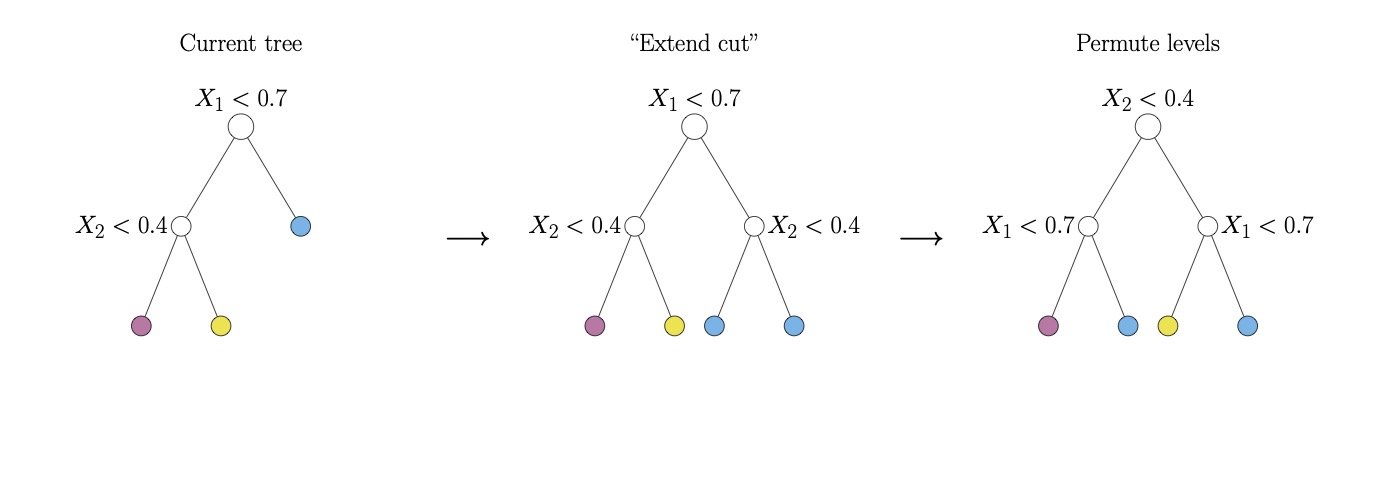}
\caption{The associated trees that correspond to the extension in Figure \ref{fig:p2_partition}. The rightmost tree corresponds to permuting the root node with all nodes in the first level.}
\label{fig:p2_process}
\end{figure}

The left and right partitions shown in Figure~\ref{fig:p2_partition} correspond, respectively, to the left and middle trees of Figure~\ref{fig:p2_process}.
Notice further that the middle tree in Figure~\ref{fig:p2_process} has the same special form as discussed above: it is a full binary tree with the same decision rule at all nodes at the same level. 
We can therefore form an equivalent tree by permuting the decision rules across levels, yielding the tree on the right of Figure~\ref{fig:p2_process}.
Unlike in the earlier examples, the left and right trees in Figure~\ref{fig:p2_process} represent different partitions of $[0,1]^{2}$: the left tree partitions the space into 3 rectangles while the right partitions it into 4.
However, because they induce the same partition of the points -- and hence, same likelihood values -- we view them as equivalent.

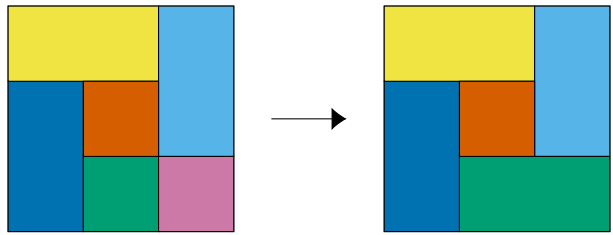
\begin{figure}[H]
\centering
\begin{tikzpicture}
\draw (0,0) rectangle (3,3);
\draw[fill = myColor6] (0,0) rectangle (1,2);
\draw[fill = myColor5] (0,2) rectangle (2,3);
\draw[fill = myColor3] (2,1) rectangle (3,3);
\draw[fill = myColor4] (1,0) rectangle (2,1);
\draw[fill = myColor7] (1,1) rectangle (2,2);
\draw[fill = myColor8] (2,0) rectangle (3,1);

\draw [->] (3.5,1.5)--(4.5,1.5);

\begin{scope}[shift={(5,0)}]
\draw (0,0) rectangle (3,3);
\draw[fill = myColor6] (0,0) rectangle (1,2);
\draw[fill = myColor5] (0,2) rectangle (2,3);
\draw[fill = myColor3] (2,1) rectangle (3,3);
\draw[fill = myColor4] (1,0) rectangle (3,1);
\draw[fill = myColor7] (1,1) rectangle (2,2);

\draw [opacity = 0] (0,0)--(0,-0.5);

\end{scope}
\end{tikzpicture}
\caption{An example of a line which cannot be contracted. The right image cannot be constructed via trees since the decision at the root must create a line spanning opposite sides, but no such spanning line exists.}
\label{fig:invalid}
\end{figure}

To implement an MCMC algorithm, moves must be reversible. Hence it is worth noting that not only can one extend lines in the box diagram, but also one can delete them. It is important when deleting a line, that the resulting diagram still corresponds to a tree. In Figure~\ref{fig:invalid}, the right diagram does not correspond to a tree. We conjecture that the condition for a diagram to correspond to a tree is that the diagram must have at least one vertical or horizontal line that goes across the whole diagram, and this same condition must hold recursively in each of the halves that this spanning line creates.

\subsection{Next steps}

In higher dimensions (i.e., $p > 2$), constructing equivalent trees is somewhat harder.
However, if the tree is a sub-tree of a full binary tree that uses the same decision at each node at a given level, then we can still permute the decision rules across levels as in Figure~\ref{fig:full_tree}.
Similarly, if we can safely extend or contract cuts as in Figure~\ref{fig:p2_partition}, then we can form equivalent trees in higher dimensions.


\section{Phase transition in 3-way contingency tables}\label{sec:transition}

This section arose from the Questions and Consulting seminar by Hanbaek Lyu and subsequent discussions amongst Miles Bakenhus, Elizabeth Gross, Max Hill, Vishesh Karwa,  and Sonja Petrovi\'c. 
The sharp phase transition phenomenon relates to two distributions on the space of tables with fixed marginal totals, also known as the fiber of a log-linear model with those marginals as sufficient statistics. In this problem, the question is how well the hypergeometric distribution can approximate the uniform distribution on the fiber. In the case of two-dimensional tables, the answer toggles between `really well' and `rather poorly' at a phase transition point that can be explained by a condition on the margins of the table. This project translates the setup for the problem into the contingency table language used in algebraic statistics, allowing us to extend the basic notions needed to generalize the results to multi-way tables. 

\subsection{Introduction}
 Log-linear models for cross-classified categorical data -- contingency tables -- have a long history in statistics  (\cite{BishopHollandFienberg}, \cite{Fienberg}) and appear in a broad variety of applications including ecology, biology, educational testing, and network science.  These models were also some of the first studied in the modern algebraic statistics literature that took off in the 1990s, leading to a vast literature on  sampling using algebraic and hybrid techniques, testing model fit using exact conditional tests,  likelihood geometry of log-linear models, the existence of maximum likelihood estimators (MLEs), and connections to classical graphical models and, more recently, colored graphical models.   While there are many flavors and variations of how contingency tables appear in the algebraic statistics literature, one typically studies models on  tables under the multinomial, Poisson, or product-multinomial sampling scheme.  It is well known that the MLEs are the same under all three sampling schemes, see for example \cite[Chapter 3]{Fienberg}. 

 In this work we consider another sampling scheme,  namely, geometric. This distribution makes a significant appearance in another vast collection of literature  on contingency tables, here referred to as `the phase-transition literature,' which is detailed in the next section. 
One common thread that appears in both  the categorical data analysis and phase transition  literature is the use of zero-margin tables, commonly referred to as \emph{moves}, to sample the space of tables given fixed marginals. The algebraic statistics literature  derives collections of such moves called \emph{Markov bases} using techniques in computational algebraic geometry and combinatorial commutative algebra. Markov bases contain moves guaranteed to connect all sets of tables for a given choice of margins. Markov bases are theoretically defined and studied for many models and types of table margins (see \cite{MB25years} for a recent overview and \cite{DobraSullivant} whose introductory section reviews the early theoretical considerations and statistical applications). 

The table margins are sufficient statistics when considering discrete exponential family models.  This space of tables  with fixed values of sufficient statistics is called \emph{the fiber of the log-linear model}. Sampling fibers provides a bona-fide algorithm for testing model goodness of fit of every such exponential family. Fibers are reference sets for the sampling, while the desired distribution for the exact conditional test is the conditional distribution of tables given the marginals. For the Poisson and multinomial sampling schemes, this distribution is hypergeometric. 
On the other hand, 
the uniform distribution on the contingency tables also plays a central role in testing for the presence of interactions in co-occurrence tables, particularly co-occurrence (between species $i$ and habitat $j$) tables arising in
ecology \citep{connor1979assembly, snijders1991enumeration, gotelli2000null, chen2005sequential}. In these contexts, the uniform distribution subject to the margin is taken as a null hypothesis of no interaction among species without any particular modeling assumption.  \cite{diaconis1985testing} use the uniform distribution on the fiber for \emph{conditional volume tests} under the  multinomial sampling scheme.   Further, we will  see  in Section~\ref{sec:3way} that if the cells follow a geometric distribution, then the log-linear model's conditional distribution given the set of margins is, in fact, uniform (see \ref{how uniform distribution arises}). 

In the next section we focus on two-dimensional tables and discuss the previous work on phase transitions related to {uniformly} sampling the space of tables given marginal totals. In Section~\ref{sec:3way} we compute a starting set of examples for three-dimensional tables, and close with a discussion Section~\ref{sec:algebraic} on connections to algebraic statistics.

	\subsection{2-way contingency tables}

	\subsubsection{Phase transition in random two-way tables with given margins}
 
	Two-way contingency tables (CTs) 
    are $m\times n$ matrices of nonnegative integer entries with prescribed row sums $\r=(a_{1},\cdots,a_{m})$ and columns sums $\cc=(b_{1},\cdots,b_{n})$ called \textit{margins}, where by $\mathcal{M}(\r,\cc)$ we denote the set of all such tables. They are fundamental objects in statistics for studying dependence structure between two or more variables and also correspond to bipartite multi-graphs with given degrees and play an important role in combinatorics and graph theory, see e.g.~\cite{barvinok2009asymptotic}. \textit{Counting} their number $|\mathcal{M}(\r,\cc)|$ and  \textit{sampling} an element from $\mathcal{M}(\r,\cc)$ uniformly at random are two fundamental problems concerning CTs with many connections and applications to other fields \citep{branden2020lower} (e.g., testing hypothesis on co-occurrence of species in ecology \citep{connor1979assembly}). A historic guiding principle to these problems is the \textit{independent heuristic}, which was introduced by I.~J.~Good as far back as in~1950 \citep{good1950probability}. The heuristic states that the constraints for the rows and columns of the table are asymptotically independent as the size of the table grows to infinity. This yields a simple yet surprisingly accurate formula that approximates the count $|\mathcal{M}(\r,\cc)|$. The independence heuristic also implies the hypergeometric (or Fisher-Yates) distirbution should approximate the uniform distribution on $\mathcal{M}(\r,\cc)$.
	
	\begin{figure*}
		\centering
		\vspace{-0.5cm}
		\includegraphics[width = 1\textwidth]{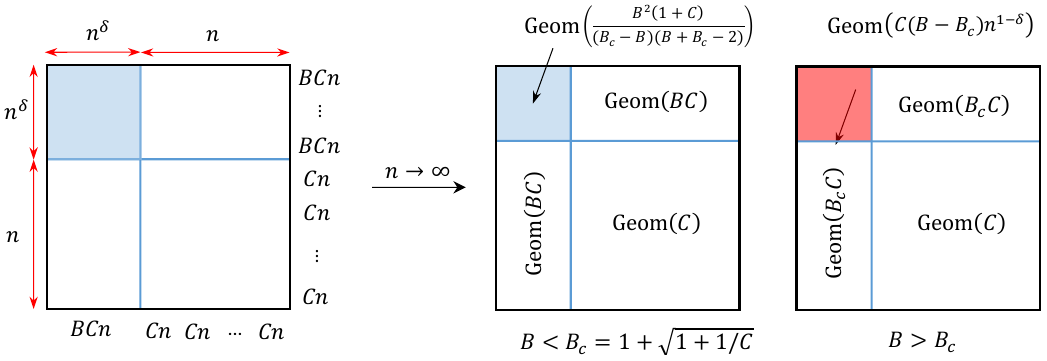}
		\caption{ (Left) Contingency table with parameters $n,\delta,B$ and $C$.
			First $\lfloor n^{\delta} \rfloor$ rows and columns have margins $\lfloor Cn \rfloor$,
			the last $n$ rows and columns have margins $\lfloor BCn \rfloor$. (Right) Limiting distributions of the entries in the uniform contingency table $X$ in the subcritical $B<B_{c}$ (left) and supercritical $B>B_{c}$ (right) regimes for thick bezels $1/2<\delta<1$. $\textup{Geom}(\lambda)$ denotes geometric distribution with mean $\lambda$. 
		}
		\label{fig:CT_thm}
	\end{figure*}

	Both of these implications of the independent heuristic have been verified when the margins are constant or have a bounded ratio close to one \citep{canfield2010asymptotic}. However, when the margins are far from being constant,  \citet{barvinok2010does} conjectured that there is a drastic difference between the uniform and hypergeometric distribution on $\mathcal{M}(\r,\cc)$, and  \citet{barvinok2012asymptotic} showed that the independent heuristic gives a large undercounting of CTs. The work of \citet{dittmer2019phase} and \citet{lyu2020number} provided the first complete answer to this puzzle; \textit{CTs exhibit a sharp phase transition when the heterogeneity of margins exceeds a certain critical threshold.} For instance, the hypergeometric distribution correctly estimates the uniform distribution for $B<B_{c}=1+\sqrt{1+1/C}$, but does so drastically differently for $B>B_{c}$ (see Figure \ref{fig:CT_thm}).  Such a sharp phase transition gives a probabilistic answer to the statistical question of why sampling a uniformly distributed CT is hard. This result settles Barvinok's conjecture \citep{barvinok2010does} (except for the special case of $\delta\in [0,1/2]$). \citeauthor{lyu2020number}  obtained a similar phase transition result for CTs, this time from the perspective of the counting problem. Roughly speaking,  the rows and columns of CTs are asymptotically independent (and hence the independence heuristic is correct) when the ratio $B$ between the two margins is strictly less than the critical threshold $B_{c}$, but suddenly they become positively correlated as soon as $B$ exceeds $B_{c}$ and the independence heuristic gives exponential undercounting.

	\subsubsection{Barvinok's typical table and the mechanism of phase transition}
	
	The key insight in \citet{dittmer2019phase} and \citet{lyu2020number} is that the uniformly random CT, say $X$, with given margins concentrates around a deterministic table called the `typical table', a notion first introduced by  \citet{barvinok2010does}. Roughly speaking, this is the $m\times n$ \textit{real-valued} table with margin $(\r,\cc)$ that maximizes a `geometric entropy function.' More precisely, let $\mathcal{P}(\r,\cc)\subseteq \R_{\ge 0}^{mn}$ denote the transportation polytope for margin $(\r,\cc)$. 
	For each $X=(X_{ij})\in \mathcal{P}(\r,\cc)$, define a strictly concave function 
	\begin{align}\label{eq:def_g}
		g(X)  =  \sum_{1\le i,j\le n}  f(X_{ij}), \quad \textup{where $f(x)  =  (x+1)\log(x+1) - x\log x$.}
	\end{align}
	The \textit{typical table} $Z\in \mathcal{P}(\r,\cc)$ for $\mathcal{M}(\r,\cc)$ is defined to be the unique maximizer of $g$ among all real-valued tables with margin $(\r,\cc)$:
	\begin{align}\label{eq:def_typical_Z}
		Z  =  \argmax_{X\in \mathcal{P}(\r,\cc)}  g(X).
	\end{align}

	The underlying mechanism of the sharp phase transition of uniformly random CTs with Barvinok margin (in Fig. \ref{fig:CT_thm}) established by 	 \citeauthor{dittmer2019phase} and \citeauthor{lyu2020number} is the sharp phase transition of the typical table $Z$, which can be shown by analyzing how the solution of the strictly concave optimization problem \eqref{eq:def_typical_Z} that defines the typical table changes as one varies the margin.

	\subsection*{Open problems} 

 Here we state three open problems related to the phase transition problem.
 
	\begin{problem}
		The phase transition in uniformly random CTs with Barvinok margin (Fig. \ref{fig:CT_thm}) is established for the `thick bezel' case $\delta>1/2$. Can one establish a similar phase transition for the `thin bezel' case $0\le \delta\le 1/2$? (c.f., Phase transition in the typical table is established for all $\delta\in [0,1]$, see \citet[Lem. 5.1]{dittmer2019phase}.)
	\end{problem}
	
	\begin{problem}\label{problem:3d_typical}
		Can we show the phase transition of the typical table when the margins assume three or more distinct values? For instance, three values  $An, Bn, Cn$ for margins for $3\times 3$ block CT of size $n\times n$, phase transition in functions of $B/A$ and $C/A$ (see \cite[Lem. 5.1]{dittmer2019phase}). In general, can one characterize all phase transitions in typical tables with respect to margin?
	\end{problem}

	\begin{problem}
		Can one show a similar sharp phase transition behavior in \cite{dittmer2019phase}  for multi-way CTs? For example, consider $n\times n\times n$ contingency tensors with margins assuming two values $BCn$ and $Cn$. \textit{Possible approach}: Develop a parallel `typical tensor' theory for uniformly random contingency tensors, and show phase transition in typical tensors as one varies the margin. Use `transference principle' to derive the behavior of uniformly random CTs from the underlying statistical model with independent entries. 
	\end{problem}

\subsection{3-way  tables with plane-sum margins}\label{sec:3way}

In this section, we provide a first study for ongoing work on sharp phase transition in 3-way contingency tables and Markov bases  \citep{Bakenhus2024sharp}. This includes the statement of a preliminary result in sharp phase transition on 3-way plane-sum contingency tables (Theorem \ref{thm:3CT_phase_transition}) and a sketch of proof. For more details we refer the interested readers for the upcoming full paper \cite{Bakenhus2024sharp}.  

For $3$-way tables, one can define more than one log-linear model. Here we consider  the \emph{model of independence}, which in log-linear model notation from \cite{Fienberg} is given by the margins $[1][2][3]$.  We will refer to these margins as $1$-margins, or plane-sums of the table.
This model is decomposable and as such is known to have a quadratic Markov basis \citep{Dobra}. 
In other words, the set of contingency tables $Y=(Y_{ijk})$ with fixed $1$-margins $Y_{i++}$, $Y_{+j+}$, and $Y_{++k}$ is \emph{connected} by moves containing exactly two $+1$'s and two $-1$'s, arranged in pairs of levels of the table so that the $1$-margins are zero.  
It is worth noting that the model of independence stands in stark contrast with the  no-three-factor interaction model \citep{Fienberg}, in which one uses $2$-margins---line sums, rather than plane sums--- as sufficient statistics. In that model, Markov bases can be arbitrarily complicated \citep{DeLoeraOnnBadNews}.

\subsubsection{The model of independence}

Fix a base measure $\mu$ on $\mathbb{Z}_{\ge 0}$. For each exponential tilting parameter $\theta\in \R$, define the exponentially tilted measure $\mu_{\theta}$ by 
\begin{align}
\frac{d\mu_{\theta}}{d\mu}(x) = e^{-\theta x -\psi(\theta)},\quad \psi(\theta):= \log\left(\sum_{k=0}^{\infty}e^{-k\theta} \mu(k) \right),
\end{align}
where $\psi$ above is the log partition function. Let $\Theta:=\{\theta\,:\, \psi(\theta)<\infty\}$, which is the set of exponential tilting parameters that makes the tilted measure $\mu_{\theta}$ a probability measure. Note that if $\theta\in \Theta$ and $X\sim \mu_{\theta}$, 
\begin{align}
    \psi'(\theta) = \E_{\mu_{\theta}}[X] \quad \textup{and} \quad  \psi''(\theta) = \textup{Var}_{\mu_{\theta}}(X)>0. 
\end{align}

The primary example is when the base measure $\mu$ is the counting measure, in which case $\Theta=(0,\infty)$  and $\mu_{\theta}$ becomes the geometric distribution on nonnegative integers with success probability $1-e^{-\theta}$, which we denote $\text{Geom}(1-e^{\theta})$: 
\begin{align}
    \P\left(\text{Geom}(1-e^{\theta})=k\right)  &= e^{k \theta}(1-e^{\theta}) \quad \textup{for $k=0,1,\dots$} \\
    & = e^{k\theta - \psi(\theta)}. 
\end{align}
Note that in this case,
\begin{align}\label{eq:psi_mean}
  \E\left[\text{Geom}(1-e^{-\theta}) \right] = \psi'(\theta) = \frac{e^{-\theta}}{1-e^{-\theta}} = \frac{1}{e^{\theta}-1}. 
\end{align}

Now define a model of $n_{1}\times n_{2}\times n_{3}$ table $Y=(Y_{ijk})$ with independent entries with marginal distribution
\begin{align}\label{eq:random_table_model1}
  Y_{ijk} \sim \mu_{\theta_{ijk}}, 
\end{align}
where $\theta_{ijk}$ is the exponential tilting parameter for the $(i,j,k)$ entry. We further make a `rank-1' assumption to the model \eqref{eq:random_table_model1}. That is, fix vector parameters $\alpha\in \R^{n_{1}}_{\ge 0}$, $\beta\in \R^{n_{1}}_{\ge 0}$, and $\gamma\in \R^{n_{1}}_{\ge 0}$. Then we assume that each exponential tilting parameter $\theta_{ijk}$ is given as the sum 
\begin{align}\label{eq:random_table_model_rank1}
    \theta_{ijk} = \alpha_{i} +\beta_{j} + \gamma_{k}. 
\end{align}
In this case, we denote 
\begin{align}
    Y=(Y_{ijk})\sim \mu_{\alpha,\beta,\gamma}. 
\end{align}

This model is  the  hierarchical log-linear specification of the model of independence \citep{BishopHollandFienberg}, whose sufficient statistics are one-dimensional table marginals.  
To see this, we note that the model can be written in the exponential family form as follows:
\begin{align}
\P\left( \{y_{ijk}\}\right) &= \prod_{ijk} e^{y_{ijk}\theta_{ijk}}\left(1-e^{\theta_{ijk}}\right)\\
&=\exp\left(\sum_{ijk}y_{ijk}\theta_{ijk} - \sum_{ijk} \psi(\theta_{ijk})\right) \\
&=\exp\left( \sum_i \alpha_i \left(\sum_{jk}y_{ijk}\right) + \sum_j \beta_j \left(\sum_{ik}y_{ijk}\right) + \sum_k \gamma_k \left(\sum_{ij}y_{ijk}\right) - \sum_{ijk} \psi(\theta_{ijk})\right) \\
&=\exp\left( \sum_i \alpha_i y_{i\bullet\bullet}  + \sum_j \beta_j y_{\bullet j \bullet} + \sum_k \gamma_k y_{ \bullet \bullet k} - \sum_{ijk} \psi(\theta_{ijk})\right).
\end{align}

The final equation is in the exponential family form $h(y)\exp\left(t(y)^T\eta(\theta) - \psi(\theta)\right)$ where $h(y) =1$ is the base measure, 
$t(y) = \left( \{y_{i \bullet \bullet}\}_{i=1}^{n_1}, \{y_{\bullet j \bullet}\}_{j=1}^{n_2}, \{y_{\bullet \bullet k} \}_{k=1}^{n_3}\right)$
is the vector of sufficient statistics, and $\eta(\theta) = (\{\alpha_i\}_{i=1}^{n_1}, \{\beta_j\}_{j=1}^{n_2}, \{\gamma_k\}_{k=1}^{n_3})$ is the vector of natural parameters. 

\subsubsection{Uniform conditional distribution on the space of tables} 
\label{how uniform distribution arises}

Fix vectors $a\in \Z^{n_{1}}_{\ge 0}$, 
$b\in \Z^{n_{2}}_{\ge 0}$, and $c\in \Z^{n_{3}}_{\ge 0}$ such that $\lVert \alpha \rVert_{1}=\lVert \beta \rVert_{1}=\lVert \gamma \rVert_{1}=:N$, where $N$ denotes the total sum. 
Define 
\begin{align}
\mathcal{T}(a,b,c) = \left\{ X\in \Z_{\ge 0}^{n_{1}\times n_{2} \times n_{3}} \,\bigg|\, \begin{matrix}
X_{i \bullet \bullet} = a_{i},\,  X_{ \bullet j \bullet} = b_{i},\, X_{ \bullet \bullet k} = c_{i} \\
\textup{for all $(i,j,k)\in [n_{1}]\times [n_{2}]\times [n_{3}]$}  
\end{matrix}\right\},
\end{align}
which is the set of all 3-way contingency tables with plane-sum margin $(a,b,c)$. 
This set is called the \emph{fiber} of the given marginal counts under the model of independence.

The probability of $Y$ conditional on satisfying the margin $(a,b,c)$ is uniform over $\mathcal{T}(a,b,c)$. To see this, note that for each $X=(x_{ijk})\in \mathcal{T}(a,b,c)$, the log likelihood of observing $X$ under the $Y$-model is 
\begin{align}
 & \log \, \P\left(Y=X  \right) - \sum_{i,j,k} \log \mu(x_{ijk}) \\
  &\qquad = (\alpha_{i}+\beta_{j}+\gamma_{k})  - \psi(\alpha_{i}+\beta_{j}+\gamma_{k} )\\
  &\qquad=\sum_{i} x_{i\bullet\bullet} \alpha_{i} + \sum_{j} x_{\bullet j \bullet} \beta_{j} + \sum_{k} x_{\bullet\bullet k}\gamma_{k} - \sum_{i,j,k} \psi(\alpha_{i}+\beta_{j}+\gamma_{k})\\  
  &\qquad= \sum_{i} a_{i} \alpha_{i} + \sum_{j} b_{i} \beta_{j} + \sum_{k} c_{k}\gamma_{k} - \sum_{i,j,k} \psi(\alpha_{i}+\beta_{j}+\gamma_{k})\\
  &\qquad=:\ell^{(a,b,c)}(\alpha,\beta,\gamma).
\end{align}
Notice that the conditional log-likelihood $\ell^{(a,b,c)}(\alpha,\beta,\gamma)$ defined above does not depend on the particular choice of $X\in \mathcal{T}(a,b,c)$, but only to the margin $(a,b,c)$ and the exponential tilting parameter $(\alpha,\beta,\gamma)$. Therefore, if the base measure $\mu$ is uniform (i.e., the counting measure), the log-likelihood is uniform over the fiber $\mathcal{T}(\r,\cc)$. Thus, in this case, the law of $Y$ conditional on being in $\mathcal{T}(a,b,c)$ is uniform. 

In fact, a similar statement is true for general $k$-way contingency table models under  the geometric sampling scheme. For completeness, we provide a short proof of this fact: 
\begin{lemma}
    Let $\bm{y}$ denote a $k$-way contingency table in its vectorized form, and let $$\P(\bm{Y}=\bm{y};\theta) = h(\bm{y}) \exp\left(t(\bm{y})^T\eta(\bm{\theta}) - \psi(\bm{\theta})\right)$$ define a exponential family model on $\bm{y}$ where $h(\bm{y})$ is the base measure, $\theta$ is a vector of parameters, $t(\bm{y})$ is the vector of sufficient statistics. Then 
\begin{align}
    \P\left(\bm{Y} = \bm{y}|t(\bm{y}) = \bm{b}\right) = \frac{h(\bm{y}) }{\sum_{\bm{y}' \in \mathcal{T}(\bm{b})}h(\bm{y}')} 
\end{align}
where $\mathcal{T}(\bm{b}) = \{\bm{y}: t(\bm{y}) = \bm{b}\}$ is the set of all tables whose sufficient statistics are equal to $\bm{b}$.

Under the Poisson and multinomial sampling schemes on the cells, this conditional distribution on the fiber is hypergeometric. Under geometric, the conditional distribution on the fiber is the uniform distribution. 
\end{lemma} 
\begin{proof}
Note that the distribution of the sufficient statistics has the following form: 
\begin{align}
    \P\left(t(\bm{y}) = \bm{b}\right) &= \sum_{\bm{y}' \in \mathcal{T}(\bm{b})}\P\left(\bm{Y} = \bm{y}'\right) \\
    &= \sum_{\bm{y}' \in \mathcal{T}(\bm{b})} h(\bm{y}') \exp\left( t(\bm{y'})^T \eta(\theta) - \psi(\theta) \right)\\
    &= \exp\left( \bm{b}^T \eta(\theta) - \psi(\theta) \right) \left(\sum_{\bm{y}' \in \mathcal{T}(\bm{b})} h(\bm{y}')\right) 
\end{align}
Hence, we have,
\begin{align}
    \P\left(\bm{y} |t(\bm{y}) = \bm{b}\right) &= \frac{\P\left(\bm{y}\right)}{\P\left(t(\bm{y})  = \bm{b}\right)} \mbox{ if } \bm{y} \in \mathcal{T}(\bm{b}), 0 \mbox{ otherwise} \\
    &=\frac{h(\bm{y}) \exp\left(\bm{b}^T\eta(\bm{\theta}) - \psi(\bm{\theta})\right)}{\exp\left( \bm{b}^T \eta(\theta) - \psi(\theta) \right) \left(\sum_{\bm{y}' \in \mathcal{T}(\bm{b})} h(\bm{y}')\right)} = \frac{h(\bm{y}) }{\sum_{\bm{y}' \in \mathcal{T}(\bm{b})}h(\bm{y}')}. 
\end{align}

It is now a direct consequence of this fact that when the underlying sampling scheme is a Poisson, or multinomial distribution, the base measure $h(\bm{y}) = \frac{N!}{\prod_{i} y_i!}$, where $N = \sum_{i} y_i$ and $i$ loops over the multi-index $\mathcal I$ of the k-way contingency table. The corresponding distribution on the fiber is hyper-geometric:
\begin{align}
\frac{\left(\prod_i y_i!\right)^{-1}}{\sum_{\bm{y}'\in \mathcal{T}(\bm{b})} (\prod_i y_i' )^{-1}}
\end{align}
When the sampling scheme is geometric, as in our case, then the base measure is $h(\bm{y}) = 1$, then the corresponding conditional distribution is uniform:
\begin{align}
\frac{1}{\sum_{\bm{y}'\in \mathcal{T}(\bm{b})} 1} = \frac{1}{|\mathcal{T}(\bm{b})|}. 
\end{align}
\end{proof}

\subsubsection{The maximum likelihood problem}

The MLE for the exponential tilting parameter $(\alpha,\beta,\gamma)$ given margin $(a,b,c)$ is given by 
\begin{align}\label{eq:MLE_problem}
    \max_{\alpha>0, \beta>0, \gamma>0} \left[ \ell^{(a,b,c)}(\alpha,\beta,\gamma) =  \sum_{i} a_{i} \alpha_{i} + \sum_{j} b_{i} \beta_{j} + \sum_{k} c_{k}\gamma_{k} - \sum_{i,j,k} \psi(\alpha_{i}+\beta_{j}+\gamma_{k}) \right]. 
\end{align}
Note that the log-likelihood function $\ell^{(a,b,c)}$ is strictly concave since 
\begin{align}
    \psi''(\theta) = \Var(\Geom(1-e^{\theta}))>0.
\end{align}
Thus if a solution for \eqref{eq:MLE_problem} exists, then it is unique and is a critical point. Thus, solving the optimization problem \eqref{eq:MLE_problem} is equivalent to solving the following MLE equations 
\begin{align}
   \begin{cases}
    &\sum_{j,k} \psi'(\alpha_{i}+\beta_{j}+\gamma_{k}) = a_{i} \quad \textup{$i=1,\dots,n_{1}$} \\
    &\sum_{k,i} \psi'(\alpha_{i}+\beta_{j}+\gamma_{k}) = b_{j}  \quad \textup{$j=1,\dots,n_{2}$} \\
    &\sum_{i,j} \psi'(\alpha_{i}+\beta_{j}+\gamma_{k}) = c_{k}   \quad \textup{$k=1,\dots,n_{3}$}. 
   \end{cases}
\end{align}
Using \eqref{eq:psi_mean}, the above is equivalent to 
\begin{align}\label{eq:MLE_eqs}
   \begin{cases}
    &\sum_{j,k} \frac{1}{\exp(\alpha_{i}+\beta_{j}+\gamma_{k})-1} = a_{i} \quad \textup{$i=1,\dots,n_{1}$} \\
    &\sum_{k,i} \frac{1}{\exp(\alpha_{i}+\beta_{j}+\gamma_{k})-1} = b_{j}  \quad \textup{$j=1,\dots,n_{2}$} \\
    &\sum_{i,j} \frac{1}{\exp(\alpha_{i}+\beta_{j}+\gamma_{k})-1} = c_{k}   \quad \textup{$k=1,\dots,n_{3}$}. 
   \end{cases}
\end{align}
More concisely, this is simply requiring that the expected table $\E[Y]$ satisfies the margin $(a,b,c)$: 
\begin{align}\label{eq:def-of-Z}
    Z:=\E[Y]\in \mathcal{T}(a,b,c). 
\end{align}

This corresponds to the standard critical equation setup in loglinear models, and here the sufficient statistics are the table plane sums. 

\subsubsection{The computational problem: plane sums of  3-way tables for 3-way Barvinok margin}

Now we consider the three-way plane-sum Barvinok margin for $(n\times n \times n)$ three-way contingency tables: 
\begin{align}\label{eq:3way_barvinok_margin}
    a = b = c = (Bn^{2},n^{2},\dots,n^{2})\in \R^{n}. 
\end{align}
The resulting MLE problem is to find the parameters $(\alpha,\beta,\gamma)$ such that the plane sums of the expected table $Z$ as in \eqref{eq:def-of-Z} match the values $Bn^{2}$ or $n^{2}$.
This amounts to satisfying the MLE equations in \eqref{eq:MLE_eqs} for the three-way Barvinok margin \eqref{eq:3way_barvinok_margin}. By symmetry, without loss of generality we can assume
\begin{equation}\label{eq:parameter-symmetry}
  \alpha = (\alpha_{1},\alpha_{2},\ldots,\alpha_{2}) = \beta = \gamma.
\end{equation}

Slices can be cut from $Z$ in each of three directions, however due to \eqref{eq:parameter-symmetry}, $Z$ is symmetric and hence it is sufficient to choose just one direction, as the plane sums of slices cut in the other two directions will be the same. Therefore it will suffice to compute the plane sums $Z_{\bullet\bullet k}$,  for $k\in [n]$. Observe that for each $k$, the slice of $Z$ corresponding to the sum $Z_{\bullet\bullet k}$ can be decomposed into four blocks of entries: a $1\times 1$ block, an $(n-1)\times(n-1)$ block, and two blocks with dimensions $1\times (n-1)$ and $(n-1)\times 1$. Due to \eqref{eq:parameter-symmetry}, all entries in a block have the same value. Using this observation, it follows that the MLE equation \eqref{eq:MLE_eqs} in this case reduces to 0
\begin{equation}\label{eq:5}
  Z_{\bullet \bullet k} = Z_{11k} + 2(n-1)Z_{12k} + (n-1)^{2}Z_{22k} = 
        \begin{cases}
        Bn^{2} \quad:\quad k=1,\\
        n^{2} \quad:\quad k \geq 2
        \end{cases}.
\end{equation}
Letting
\begin{equation*}
  P:=e^{\alpha_{1}} \quad \text{and} \quad Q:=e^{\alpha_{2}},
\end{equation*}
we have 
\begin{align}\label{eq:Z_PQ}
Z_{111} = \frac{1}{P^{3}-1}, \quad Z_{121} =Z_{112}= \frac{1}{P^{2}Q-1}, \quad Z_{221} =Z_{122} = \frac{1}{Q^{2}P-1}, \quad Z_{222}=\frac{1}{Q^{3}-1}. 
\end{align}
Since $Z_{ijk}\ge 0$, we must have $P,Q\ge 1$. Also, this change of variables rewrites  \eqref{eq:5} as
\begin{align}\label{eq:barvinok_MLE_PQ}
\begin{cases}
  \frac{1}{P^{3}-1}+2(n-1)\frac{1}{P^{2}Q - 1} + \frac{(n-1)^{2}}{Q^{2}P - 1} = Bn^{2} \\
  \frac{1}{P^{2}Q-1}+2(n-1)\frac{1}{PQ^{2} - 1} + \frac{(n-1)^{2}}{Q^{3} - 1} = n^{2}.
 \end{cases}
\end{align}
Since all terms involved are nonnegative, dropping the first two terms in the second equation in \eqref{eq:barvinok_MLE_PQ}, it follows that $\frac{(n-1)^{3}}{Q^{3}-1}\le n^{3}$, so 
\begin{align}
    Q \ge \left(\left( \frac{n}{n-1} \right)^{3} +1\right)^{1/3} =2+O(n^{-1}). 
\end{align}
Thus $Q^{3}>3/2$ for all sufficiently large $n$. It follows that $Z_{ijk}=O(1)$ whenever $(i,j,k)\ne (1,1,1)$. Reusing the second equation in \eqref{eq:barvinok_MLE_PQ} with this fact, we have 
\begin{align}\label{eq:Q_asymptotics}
    \frac{(n-1)^{2}}{Q^{3}-1} = n^{2} - O(n), 
\end{align}
so we deduce 
\begin{align}
    Q^{3} = 2 + O(n^{-1}). 
\end{align}
Thus, for all values of the ratio parameter $B$, $Q\sim 2^{1/3}$. On the contrary, as we will see shortly, the asymptotic of $P$ will depend on $B$ critically. 

\subsubsection{Sharp phase transition in the maximum-likelihood expected table}

Below is the result on sharp phase transition in the expected table under the MLE model with 3D Barvinok margin $(Bn^{2},n^{2},\dots, n^{2})$. This result is a three-way extension of \cite[Lem. 5.1]{dittmer2019phase}, which partially addresses Problem \ref{problem:3d_typical}.  

\begin{theorem}\label{thm:3CT_phase_transition}
    Let $Z=Z^{\r,\cc}=\E_{\alpha,\beta,\gamma}[Y]$, where $Y\sim  \mu_{\alpha,\beta,\gamma}$ and $(\alpha,\beta,\gamma)$ is an MLE for the rank-1 model for the plane-sum margin $(a,b,c)$ with $a=b=c=(Bn^{2}, n^{2},\dots,n^{2})\in \R^{n}$. For $B_{c}:=\frac{1}{2^{2/3}-1}$, the followings hold:
    \begin{description}
        \item[(i)] (Subcritial regime) Suppose $B<B_{c}$. Then 
        \begin{align}
            &Z_{111} = \frac{1}{\left(\frac{B^{-1}+1}{B_{c}^{-1}+1} \right)^{3}-1} + O(n^{-1}), \qquad Z_{121}=Z_{112} = \frac{1}{2^{1/3}\left(\frac{B^{-1}+1}{B_{c}^{-1}+1} \right)^{2}-1} + O(n^{-1}), \\
            &Z_{221}=Z_{122} = \frac{1}{2^{2/3}\left(\frac{B^{-1}+1}{B_{c}^{-1}+1} \right)-1} + O(n^{-1}),\qquad Z_{222}  = 1 + O(n^{-1}). 
        \end{align}

        \item[(ii)] (Supercritial regime) Suppose $B>B_{c}$. Then 
        \begin{align}
           & Z_{111} = (B-B_{c})n^{2} - O(n), \qquad Z_{121}=Z_{112} = \frac{1}{2^{1/3}-1} + O(n^{-1}), \\
           & Z_{221}=Z_{122} = \frac{1}{2^{2/3}-1} + O(n^{-1}), \qquad Z_{222}  = 1 + O(n^{-1}). 
        \end{align}
    \end{description}
\end{theorem}

\begin{proof}[Sketch of proof]
First assume $B<B_{c}$. Dropping the first two terms in the first equation in \eqref{eq:barvinok_MLE_PQ} and using the asymptotic of $Q$ in \eqref{eq:Q_asymptotics}, 
\begin{align}
\frac{(n-1)^{2}}{Q^{2}P - 1} \le Bn^{2}. 
\end{align}
Solving for $P$ and using $B<B_{c}$, we get 
\begin{align}
 P\ge  Q^{-2} \left( B^{-1}\left( \frac{n-1}{n} \right)^{2} +1 \right) \sim 2^{-2/3}(B^{-1}+1) = \frac{B^{-1}+1}{B_{c}^{-1}+1}>1.   
\end{align}
It follows that $Z_{111}=\frac{1}{P^{3}-1}=O(1)$. Hence the first two terms in the first equation in \eqref{eq:barvinok_MLE_PQ} is of order $O(n)$, so we get 
\begin{align}
\frac{(n-1)^{2}}{Q^{2}P - 1} = Bn^{2} - O(n). 
\end{align}
Consequently, we deduce the following asymptotic of $P$ in the subcritical regime:
\begin{align}
    P = \frac{B^{-1}+1}{B_{c}^{-1}+1} + O(n^{-1}). 
\end{align}
Now using \eqref{eq:Z_PQ} and the asymptotic of $Q$ in \eqref{eq:Q_asymptotics}, we conclude the asymptotics of the entries of $Z$ stated in \textbf{(i)}.

Next, assume $B>B_{c}$. Recall that taking $P\searrow 1$, the largest $n^{2}$ asymptotic that the third term in the first equation in \eqref{eq:barvinok_MLE_PQ} is at most $B_{c}n^{2}$. Also recall that the second term in the first equation in \eqref{eq:barvinok_MLE_PQ} is of order $O(n)$. Therefore, since $B>B_{c}$, it follows that the first term $\frac{1}{P^{3}-1}$ must have contributions in the $(B-B_{c})n^{2}$.  We can obtain more precise asymptotics as follows. Using the first equation in \eqref{eq:barvinok_MLE_PQ} and that $Z_{121}=O(1)$, 
\begin{align}
    \frac{1}{P^{3}-1}+ \frac{(n-1)^{2}}{Q^{2}P - 1} = Bn^{2} - O(n). 
\end{align}
Since $Q^{3}=2+O(n^{-1})$, $P^{3}=1+O(n^{-2})$, so 
\begin{align}
 Z_{111}= \frac{1}{P^{3}-1} = (B-B_{c})n^{2} - O(n).
\end{align}
From this, we get 
\begin{align}
    P^{3}=1+\frac{n^{-2}}{B-B_{c}} + O(n^{-3}). 
\end{align}
Then using \eqref{eq:Z_PQ} and the asymptotic of $Q$ in \eqref{eq:Q_asymptotics}, we can also conclude the asymptotics of the entries of $Z$ stated in \textbf{(ii)}. 
\end{proof}

\subsubsection{Algebraic statistics and phase transitions}\label{sec:algebraic}

We propose to connect the phase transition phenomena to Markov bases for sampling fibers and also to the geometry of marginal polytopes. 

First, phase transitions in two-way tables happen for so-called tame fibers, that is, 
two $m\times n$ margins that are `$\delta$-tame' (MLEs are uniformly bounded by $\delta^{-1}$). From the point of view of discrete exponential families, this means that the MLE is away from the boundary of the marginal polytope, which is defined as the convex hull of all possible sufficient statistics for the exponential family. 
This implies that the phase transition can be also interpreted on an infinite sequence of polytopes: for small table sizes, the MLE is well in the relative interior of the polytope, but as the size grows, there is a threshold when the MLE comes too close (within $\delta$) to the face of the polytope.  It would be of great interest to perform some initial simulations in this direction.

Second, the moves used to connect the fiber in the $2$-way table case provide a basis for a modified Markov chain to sample from the fiber uniformly, until the phase transition occurs with the mass blowup on the corner $(1,1)$ cell of the cube. Naturally, we expect that this phenomenon  generalizes to the $3$-way independence model when the mass concentrates on the $(1,1,1)$ cell.  
To this end, 
we close with an explicit description of Markov basis moves for the $3$-way model of independence. These are the moves analogous to the moves in the $2$-way table model that can be used to explore the fiber uniformly, under the `tame' regime before the phase transition occurs. 

The Markov moves $m=(m_{ijk})$ for the independence model on the $n_1 \times n_2 \times n_3$ table $Y = (Y_{ijk})$ have entries $m_{i_1j_1k_1} = 1$ and $m_{i_2j_2k_2} = 1$ with one of 
\[
\begin{cases} 
m_{i_1j_2k_2}=-1,\\
 m_{i_2j_1k_2}=-1,\\
\end{cases},\, \begin{cases}
m_{i_1j_2k_2}=-1,\\
m_{i_2j_2k_1}=-1,\\
\end{cases},\, \text{ or }\begin{cases} 
m_{i_2j_1k_2}=-1,\\
m_{i_2j_2k_1}=-1,\\
\end{cases}
\]
where all other entries are zero, for each $i_1,i_2\in [n_1]$, $j_1,j_2\in [n_2]$, and $k_1,k_2\in [n_3]$, where $i_1\ne i_2$, $j_1 \ne j_2$, and $k_1\ne k_2$. There are a total of 
\(\frac{3}{2}n_1n_2n_3(n_1-1)(n_2-1)(n_3-1)\)
of these moves. Fixing $i_1j_1k_1 = (1,1,1)$ leaves 
\(3(n_1 - 1)(n_2 - 1)(n_3 - 1)\)
moves that apply to the entry $Y_{111}$ in the table. Then, the ratio of applicable moves to non-applicable moves is 
\[\frac{2}{n_1n_2n_3}\to 0\, \text{ as } n_1,n_2,n_3 \to \infty.\]

Simply the number of available moves does not reveal directly any phase transition phenomena, but this is likely the case because Markov moves connect all possible fibers. In forthcoming work we will study what moves are necessary and/or inapplicable for the  particular margins are of the form $(Bn^2, n^2, n^2)$.


\bibliographystyle{abbrvnat}
\bibliography{flex_bart}  

\begin{thebibliography}{38}
\providecommand{\natexlab}[1]{#1}
\providecommand{\url}[1]{\texttt{#1}}
\expandafter\ifx\csname urlstyle\endcsname\relax
  \providecommand{\doi}[1]{doi: #1}\else
  \providecommand{\doi}{doi: \begingroup \urlstyle{rm}\Url}\fi

\bibitem[Alexandr and Heaton(2021)]{alexandr-heaton}
Y.~Alexandr and A.~Heaton.
\newblock Logarithmic {V}oronoi cells.
\newblock \emph{Algebraic Statistics}, 12\penalty0 (1):\penalty0 75--95, 2021.
\newblock ISSN 2693-2997.
\newblock \doi{10.2140/astat.2021.12.75}.
\newblock URL \url{http://dx.doi.org/10.2140/astat.2021.12.75}.

\bibitem[Allman et~al.(2015)Allman, Rhodes, Sturmfels, and
  Zwiernik]{allman2015tensors}
E.~S. Allman, J.~A. Rhodes, B.~Sturmfels, and P.~Zwiernik.
\newblock Tensors of nonnegative rank two.
\newblock \emph{Linear Algebra and Its Applications}, 473:\penalty0 37--53,
  2015.

\bibitem[Allman et~al.(2019)Allman, Ba{\~n}os, Evans, Ho{\c{s}}ten, Kubjas,
  Lemke, Rhodes, and Zwiernik]{allman2019maximum}
E.~S. Allman, H.~Ba{\~n}os, R.~Evans, S.~Ho{\c{s}}ten, K.~Kubjas, D.~Lemke,
  J.~A. Rhodes, and P.~Zwiernik.
\newblock Maximum likelihood estimation of the latent class model through model
  boundary decomposition.
\newblock \emph{Journal of Algebraic Statistics}, 10\penalty0 (1), 2019.

\bibitem[Almendra-Hern\'andez et~al.(2023)Almendra-Hern\'andez, De~Loera, and
  Petrovi\'c]{MB25years}
F.~Almendra-Hern\'andez, J.~A. De~Loera, and S.~Petrovi\'c.
\newblock Markov bases: a 25 year update.
\newblock \emph{Journal of the American Statistical Association}, 2023.

\bibitem[Bakenhus et~al.(2024)Bakenhus, Gross, Hill, Karwa, Lyu, and
  Petrovi\'c]{Bakenhus2024sharp}
M.~Bakenhus, E.~Gross, M.~Hill, V.~Karwa, H.~Lyu, and S.~Petrovi\'c.
\newblock Sharp phase transition in 3-way contingency tables and {M}arkov
  bases.
\newblock \emph{In preparation}, 2024.

\bibitem[Barvinok(2009)]{barvinok2009asymptotic}
A.~Barvinok.
\newblock Asymptotic estimates for the number of contingency tables, integer
  flows, and volumes of transportation polytopes.
\newblock \emph{International Mathematics Research Notices}, 2009\penalty0
  (2):\penalty0 348--385, 2009.

\bibitem[Barvinok(2010)]{barvinok2010does}
A.~Barvinok.
\newblock What does a random contingency table look like?
\newblock \emph{Combinatorics, Probability and Computing}, 19\penalty0
  (4):\penalty0 517--539, 2010.

\bibitem[Barvinok and Hartigan(2012)]{barvinok2012asymptotic}
A.~Barvinok and J.~Hartigan.
\newblock An asymptotic formula for the number of non-negative integer matrices
  with prescribed row and column sums.
\newblock \emph{Transactions of the American Mathematical Society},
  364\penalty0 (8):\penalty0 4323--4368, 2012.

\bibitem[Bishop et~al.(2007)Bishop, Fienberg, and
  Holland]{BishopHollandFienberg}
Y.~M. Bishop, S.~E. Fienberg, and P.~W. Holland.
\newblock \emph{Discrete multivariate analysis: Theory and practice}.
\newblock Springer Science \& Business Media, 2007.

\bibitem[Br{\"a}nd{\'e}n et~al.(2020)Br{\"a}nd{\'e}n, Leake, and
  Pak]{branden2020lower}
P.~Br{\"a}nd{\'e}n, J.~Leake, and I.~Pak.
\newblock Lower bounds for contingency tables via lorentzian polynomials.
\newblock \emph{arXiv preprint arXiv:2008.05907}, 2020.

\bibitem[Canfield and McKay(2010)]{canfield2010asymptotic}
E.~R. Canfield and B.~D. McKay.
\newblock Asymptotic enumeration of integer matrices with large equal row and
  column sums.
\newblock \emph{Combinatorica}, 30\penalty0 (6):\penalty0 655, 2010.

\bibitem[Chen et~al.(2005)Chen, Diaconis, Holmes, and Liu]{chen2005sequential}
Y.~Chen, P.~Diaconis, S.~P. Holmes, and J.~S. Liu.
\newblock Sequential {M}onte {C}arlo methods for statistical analysis of
  tables.
\newblock \emph{Journal of the American Statistical Association}, 100\penalty0
  (469):\penalty0 109--120, 2005.

\bibitem[Chipman et~al.(2010)Chipman, George, and McCulloch]{Chipman2010}
H.~A. Chipman, E.~I. George, and R.~E. McCulloch.
\newblock {BART}: {Bayesian} additive regression trees.
\newblock \emph{Annals of Applied Statistics}, 4\penalty0 (1):\penalty0
  266--298, 2010.

\bibitem[Connor and Simberloff(1979)]{connor1979assembly}
E.~F. Connor and D.~Simberloff.
\newblock The assembly of species communities: chance or competition?
\newblock \emph{Ecology}, 60\penalty0 (6):\penalty0 1132--1140, 1979.

\bibitem[De~Loera and Onn(2006)]{DeLoeraOnnBadNews}
J.~A. De~Loera and S.~Onn.
\newblock Markov bases of three-way tables are arbitrarily complicated.
\newblock \emph{Journal of Symbolic Computation}, 41:\penalty0 173--181, 2006.

\bibitem[Diaconis and Efron(1985)]{diaconis1985testing}
P.~Diaconis and B.~Efron.
\newblock Testing for independence in a two-way table: new interpretations of
  the chi-square statistic.
\newblock \emph{The Annals of Statistics}, pages 845--874, 1985.

\bibitem[Dittmer et~al.(2020)Dittmer, Lyu, and Pak]{dittmer2019phase}
S.~Dittmer, H.~Lyu, and I.~Pak.
\newblock Phase transition in random contingency tables with non-uniform
  margins.
\newblock \emph{Transactions of the AMS.}, 373:\penalty0 8313--8338, 2020.

\bibitem[Dobra(2003)]{Dobra}
A.~Dobra.
\newblock Markov bases for decomposable graphical models.
\newblock \emph{Bernoulli}, 9\penalty0 (6):\penalty0 1093--1108, 2003.

\bibitem[Dobra and Sullivant(2004)]{DobraSullivant}
A.~Dobra and S.~Sullivant.
\newblock A divide-and-conquer algorithm for generating {M}arkov bases of
  multi-way tables.
\newblock \emph{Computational Statistics}, 19:\penalty0 347--366, 2004.

\bibitem[Drton et~al.(2008)Drton, Sturmfels, and Sullivant]{drton2008lectures}
M.~Drton, B.~Sturmfels, and S.~Sullivant.
\newblock \emph{Lectures on algebraic statistics}, volume~39.
\newblock Springer Science \& Business Media, 2008.

\bibitem[Fienberg(2007)]{Fienberg}
S.~E. Fienberg.
\newblock \emph{The analysis of cross-classified categorical data}.
\newblock Springer Science \& Business Media, 2007.

\bibitem[Fienberg et~al.(2009)Fienberg, Hersh, Rinaldo, and Zhou]{fienberg2009}
S.~E. Fienberg, P.~Hersh, A.~Rinaldo, and Y.~Zhou.
\newblock \emph{Maximum likelihood estimation in latent class models for
  contingency table data}, pages 27--62.
\newblock Cambridge University Press, 2009.

\bibitem[Good(1950)]{good1950probability}
I.~J. Good.
\newblock \emph{Probability and the Weighing of Evidence}.
\newblock C. Griffin London, 1950.

\bibitem[Gotelli(2000)]{gotelli2000null}
N.~J. Gotelli.
\newblock Null model analysis of species co-occurrence patterns.
\newblock \emph{Ecology}, 81\penalty0 (9):\penalty0 2606--2621, 2000.

\bibitem[Gu(2022)]{gu2022blessing}
Y.~Gu.
\newblock Blessing of dependence: Identifiability and geometry of discrete
  models with multiple binary latent variables.
\newblock \emph{arXiv preprint arXiv:2203.04403}, 2022.

\bibitem[Gu(2023)]{gu2023generic}
Y.~Gu.
\newblock Generic identifiability of the {DINA} model and blessing of latent
  dependence.
\newblock \emph{Psychometrika}, 88\penalty0 (1):\penalty0 117--131, 2023.

\bibitem[Gu and Dunson(2023)]{gu2023bp}
Y.~Gu and D.~B. Dunson.
\newblock {Bayesian Pyramids: identifiable multilayer discrete latent structure
  models for discrete data}.
\newblock \emph{Journal of the Royal Statistical Society Series B: Statistical
  Methodology}, 85\penalty0 (2):\penalty0 399--426, 2023.

\bibitem[Hauenstein et~al.(2014)Hauenstein, Rodriguez, and
  Sturmfels]{HRS2014-mldegree-rank-constraints}
J.~Hauenstein, J.~I. Rodriguez, and B.~Sturmfels.
\newblock Maximum likelihood for matrices with rank constraints.
\newblock \emph{J. Algebr. Stat.}, 5\penalty0 (1):\penalty0 18--38, 2014.
\newblock ISSN 1309-3452.
\newblock \doi{10.18409/jas.v5i1.23}.
\newblock URL \url{https://doi.org/10.18409/jas.v5i1.23}.

\bibitem[Kim and Ro\v{c}kov\'{a}(2023)]{KimRockova2023_mixing}
J.~Kim and V.~Ro\v{c}kov\'{a}.
\newblock On mixing rates for {Bayesian CART}.
\newblock \textit{arXiv preprint 2306.00126}, 2023.

\bibitem[Kubjas et~al.(2015)Kubjas, Robeva, and Sturmfels]{KRS2015-EM}
K.~Kubjas, E.~Robeva, and B.~Sturmfels.
\newblock Fixed points {EM} algorithm and nonnegative rank boundaries.
\newblock \emph{Ann. Statist.}, 43\penalty0 (1):\penalty0 422--461, 2015.
\newblock ISSN 0090-5364,2168-8966.
\newblock \doi{10.1214/14-AOS1282}.
\newblock URL \url{https://doi.org/10.1214/14-AOS1282}.

\bibitem[Landsberg(2012)]{landsberg2012tensors}
J.~M. Landsberg.
\newblock \emph{Tensors: geometry and applications}, volume 128 of
  \emph{Graduate Studies in Mathematics}.
\newblock American Mathematical Society, Providence, RI, 2012.
\newblock ISBN 978-0-8218-6907-9.
\newblock \doi{10.1090/gsm/128}.
\newblock URL \url{https://doi.org/10.1090/gsm/128}.

\bibitem[Lyu and Pak(2020)]{lyu2020number}
H.~Lyu and I.~Pak.
\newblock On the number of contingency tables and the independence heuristic.
\newblock \emph{arXiv preprint arXiv:2009.10810}, 2020.

\bibitem[Raicu(2013)]{raicu}
C.~Raicu.
\newblock $3\times 3$ minors of catalecticants.
\newblock \emph{Mathematical Research Letters 20}, \penalty0 (745-756), 2013.

\bibitem[Rodriguez and Wang(2017)]{RW2017-mixture-independence}
J.~I. Rodriguez and B.~Wang.
\newblock The maximum likelihood degree of mixtures of independence models.
\newblock \emph{SIAM J. Appl. Algebra Geom.}, 1\penalty0 (1):\penalty0
  484--506, 2017.
\newblock ISSN 2470-6566.
\newblock \doi{10.1137/16M1088843}.
\newblock URL \url{https://doi.org/10.1137/16M1088843}.

\bibitem[Ronen et~al.(2022)Ronen, Saarinen, Tan, Duncan, and
  Yu]{Ronen2022_mixing}
O.~Ronen, T.~Saarinen, Y.~S. Tan, J.~Duncan, and B.~Yu.
\newblock A mixing time lower bound for a simplified version of {BART}.
\newblock \textit{arXiv preprint 2210.09352v1}, 2022.

\bibitem[Seigal and Mont{\'u}far(2018)]{seigal2018mixtures}
A.~Seigal and G.~Mont{\'u}far.
\newblock Mixtures and products in two graphical models.
\newblock \emph{Journal of Algebraic Statistics}, 9\penalty0 (1), 2018.

\bibitem[Snijders(1991)]{snijders1991enumeration}
T.~A. Snijders.
\newblock Enumeration and simulation methods for 0--1 matrices with given
  marginals.
\newblock \emph{Psychometrika}, 56:\penalty0 397--417, 1991.

\bibitem[Sullivant(2018)]{Sullivant-book}
S.~Sullivant.
\newblock \emph{Algebraic statistics}, volume 194 of \emph{Graduate Studies in
  Mathematics}.
\newblock American Mathematical Society, Providence, RI, 2018.
\newblock ISBN 978-1-4704-3517-2.
\newblock \doi{10.1090/gsm/194}.
\newblock URL \url{https://doi.org/10.1090/gsm/194}.

\end{thebibliography}

\end{document}